\documentclass[fleqn]{article}

\usepackage{wrapfig}
\usepackage[justification=centering]{caption}
\usepackage[top=1in, bottom=1.25in, left=1in, right=1in]{geometry}
\usepackage{fleqn}

\usepackage{bm,upgreek}
\usepackage{etoolbox}
\usepackage{xstring}
\newcommand{\gb}[1]{ 
	\StrGobbleLeft{\detokenize{#1}}{1}[\chrcodet] 
	\StrGobbleRight{\chrcodet}{1}[\chrcode] 
	\StrLeft{\chrcode}{1}[\chrfirst] 
	\IfSubStr{ABCDEFGHIJKLMNOPQRSTUVWXYZ}{\chrfirst} 
	{\boldsymbol{#1}} 
	{\boldsymbol{\csname up\chrcode\endcsname}} 
}
\newtheorem{theorem}{Theorem}[section]

\newtheorem{proposition}[theorem]{Proposition}

\newtheorem{remark}[theorem]{Remark}
\newtheorem{definition}[theorem]{Definition}

\newenvironment{proof}[1][Proof]{\begin{trivlist}
\item[\hskip \labelsep {\bfseries #1}]}{\end{trivlist}}

\newcommand{\qed}{\nobreak \ifvmode \relax \else
      \ifdim\lastskip<1.5em \hskip-\lastskip
      \hskip1.5em plus0em minus0.5em \fi \nobreak
      \vrule height0.75em width0.5em depth0.25em\fi}

\usepackage[english]{babel}
\usepackage[utf8]{inputenc}

\usepackage{titlesec}
\usepackage{titletoc}
\usepackage{lmodern}
\usepackage[T1]{fontenc}
\usepackage{textcomp}
\usepackage{amsmath}
\usepackage{amsfonts}
\usepackage{amssymb}
\usepackage{graphicx} 
\DeclareGraphicsExtensions{.pdf,.jpeg,.png}
\usepackage{caption}
\usepackage{epstopdf}
\usepackage[colorinlistoftodos]{todonotes}

\usepackage{listings}
\usepackage{color} 
\definecolor{mygreen}{RGB}{28,172,0} 
\definecolor{mylilas}{RGB}{170,55,241}
\usepackage{cite}
\usepackage{url}
\usepackage{pdfpages} 
\usepackage[colorlinks=true]{hyperref}
\hypersetup{colorlinks,
citecolor=blue,
filecolor=blue,
linkcolor=blue,
urlcolor=blue}

\usepackage{xcolor}

\usepackage{subfig}
\usepackage{longtable}
\usepackage{soul}
\usepackage{afterpage}

\title{\bf Observability of Sensorless Electric Drives}
\author{Mohamad Koteich$^{1,2}$, Gilles Duc$^2$, Abdelmalek Maloum$^1$, Guillaume Sandou$^2$\\
$^1$ Renault S.A.S., Technocentre, 78280, Guyancourt, France\\
Email: \{name.surname\}@renault.com\\
$^2$ L2S, CentraleSupelec - CNRS - Universit\'{e} Paris-Saclay, 91190, Gif-sur-Yvette, France\\
Email: \{name.surname\}@centralesupelec.fr}

\begin{document}
\setlength{\parindent}{5mm}

\maketitle

\begin{abstract}
Electric drives control without shaft sensors has been an active research topic for almost three decades. It consists of estimating the rotor speed and/or position from the currents and voltages measurement. This paper deals with the observability conditions of electric drives in view of sensorless control. The models of such systems are strongly nonlinear. For this reason, a local observability approach is applied to analyze the deteriorated performance of sensorless drives in some operating conditions. The validity of the observability conditions is confirmed by numerical simulations and experimental data, using an extended Kalman filter as observer.
\end{abstract}


\section{Introduction}
Electrical machines (EMs) have been the driving horse of the contemporary industries for more than a century \cite{tesla_aiee_88}. They found their application mainly in electric power generation, transport industry, and motor-driven systems. More recently, electric drives have become a serious competitor with combustion engines, and they are gaining more interest as Eco-friendly technologies are being sought worldwide.

\begin{figure}[!htb]
	\centering
	\includegraphics[scale = 1.2]{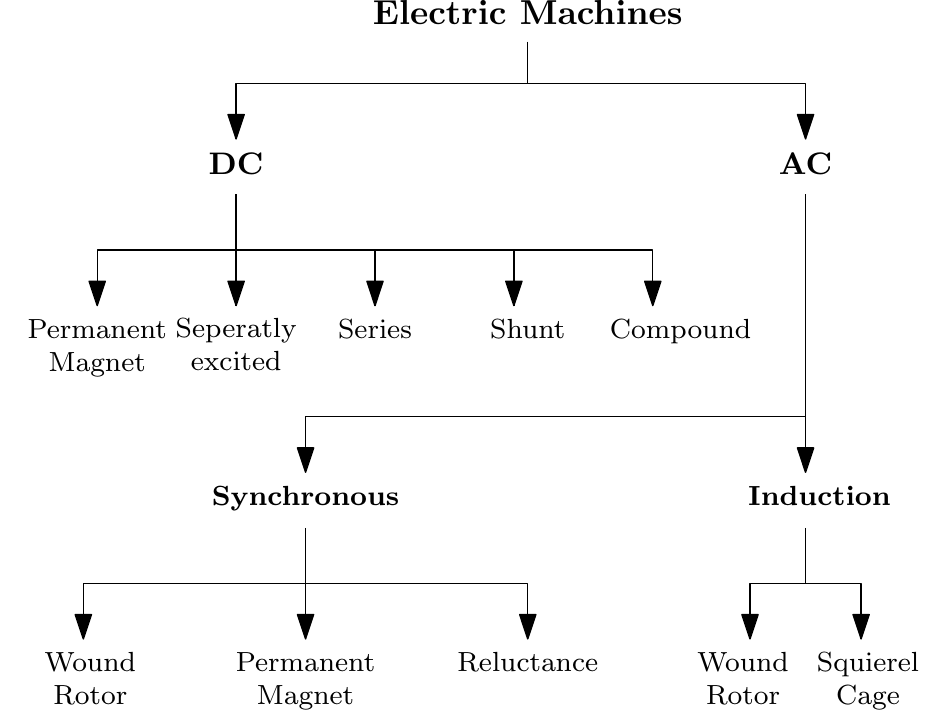}
	\caption{Electrical machines family tree}
	\label{EMs}
\end{figure}

At the beginning of the XX$^{th}$ century, most of today's electrical machines had been already invented. There exist broadly two EM families: direct-current machines (DCMs) and alternating current machines (ACMs). Figure \ref{EMs} shows a family tree of the most commonly used EMs. It shows 5 types of DCM, and 2 sub-families of ACM: synchronous machine (SM), where the rotor angular velocity is equal to the rotating magnetic field one, and induction machine (IM), where the aforementioned angular velocities are slightly different.

DC machines had dominated the industries of variable speed drives (VSD) for a while, due to their ease of control. The use of AC machines in VSD has gradually emerged due to several theoretical contributions, such as the two-reaction theory introduced by Park \cite{park_aiee_29, stanley_aiee_38} and advances in control theory, with some technological innovations, such as advances in power semiconductors and converters. Consequently, new trends of AC VSD control techniques have appeared \cite{leonhard}, such as, inter alia, the field oriented control (FOC) \cite{blaschke_siemens_72} and the direct torque control (DTC) \cite{Takahashi_tia_89}. 

\begin{figure}[!htb]
	\centering
	\includegraphics[scale = 0.9]{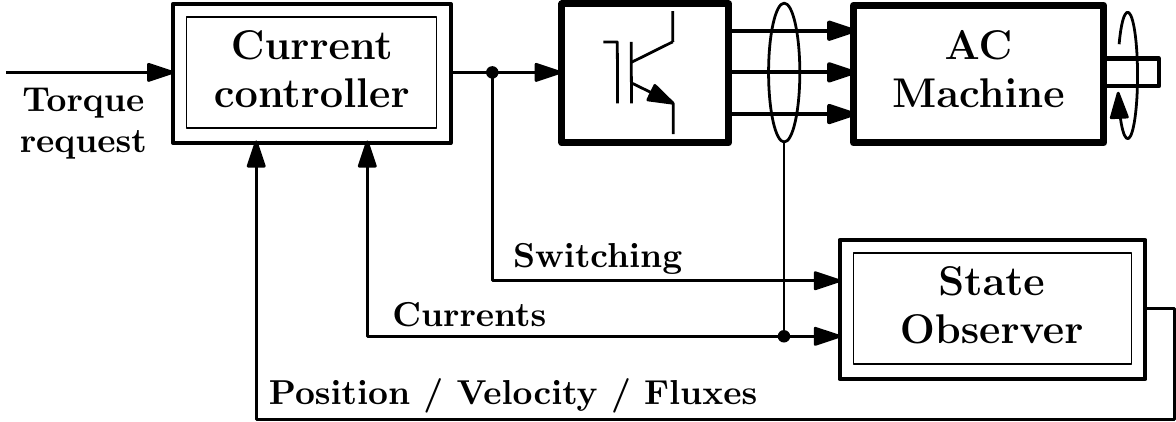}
	\caption{Schematic of the general structure of sensorless control}
	\label{sans_capteur}
\end{figure}

These new high-performance control techniques require, often, an accurate knowledge of the rotor angular speed and position. These mechanical variables are traditionally
measured using sensors, which yields additional cost of sensors and their installation and reduces the system's reliability and robustness. For these reasons, sensorless\footnote{Some engineers and researchers argue about the term \emph{sensorless}, and propose other (reasonable) terms such as \emph{self-sensing}, or \emph{mechanical-sensorless}, or \emph{encoder-less} etc. In this paper, the term \emph{sensorless} is adopted, since it has been used in some of the most important references in the field, e.g. \cite{vas, holtz_tie_06}.} control techniques have gained increasing attention recently \cite{vas}: they consist of sensing the drive's currents and voltages, and using them as inputs to an estimation algorithm (software sensor), such as a state-observer, that estimates the desired variables (figure \ref{sans_capteur}). 

After almost three decades of research and development in sensorless control, the literature today is rich in a wide variety of estimation techniques \cite{holtz_tie_06, acarnley_tie_06, finch_tie_08, pacas_iemag_11}: there exist, for instance, model-reference-adaptive based \cite{schauder_tia_92}, Kalman filter based \cite{bolognani_tia_03, barut_tie_07}, flux observer based \cite{boldea_tec_09,koteich_ecmsm_13}, back-electromotive force (EMF) observer based \cite{chen_iac_00}, sliding-mode based \cite{delpoux_tcon_14, hamida_tcon_14}, interconnected-adaptive-observer based \cite{traore_tcon_09, ezzat_tcon_11},  and signal injection based \cite{jansen_tia_95, jebai_tcon_15} estimation techniques.


\subsection{State-of-the-art of electric drives observability}

\begin{figure}[!htb]
\centering
\includegraphics[scale = 1]{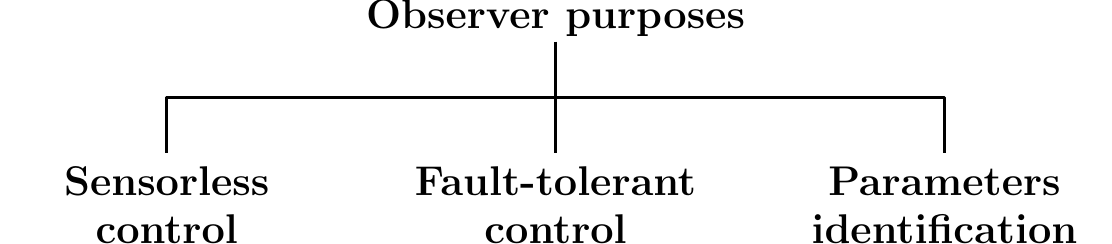}
\caption{State-observer purposes}
\label{objectifs}
\end{figure}

State-observers have been widely applied in electric drives control. Originally, the purpose was to estimate unmeasurable states, such as rotor fluxes and currents in an induction machine. However, nowadays, state observers are applied for different purposes (figure \ref{objectifs}): mainly, for replacing some sensors like in sensorless control, for sensors diagnosis and fault-tolerant control, and for parameter identification.

In this paper, only observer design in view of sensorless control is considered. Nevertheless, the results are also valid for diagnosis and fault tolerant control. For studies on the observability of some machines parameters (called sometimes \emph{identifiability}), the reader is referred to \cite{marino_10, vaclavek_tie_13} for induction drives, and to \cite{hamida_tcon_14, glumineau_15} for synchronous drives.

One limitation of the use of observer-based sensorless techniques is the deteriorated performance in some operating conditions: namely the zero- (and low-)speed operation in the case of synchronous machines (SMs), and the low-stator-frequency in the case of induction machines (IMs). Usually, this problem is viewed as a stability problem, and sometimes is treated based on experimental results. However, the real problem lies in the \emph{observability conditions} of the drive \cite{holtz_ieeeproc_02}.

AC machines are nonlinear systems, and their observability analysis is not an easy task. To the best of the authors knowledge, the first paper where the observability of an AC machine was studied in view of sensorless control is the paper by Holtz \cite{holtz_iecon_93}. In this paper, the study is based on the signal flow analysis of the machine, and the observability of the IM is viewed as a limitation on the speed estimation at low input frequencies. In 2000, the observability of the IM was studied in \cite{canudas_cdc_00} based on the local weak observability theory \cite{hermann_tac_77} for the first time. As a certain rank criterion is to be verified in this approach, the observability conditions can be formulated with analytical equations that can be examined in real time, which turns out to be extremely useful for practical implementation. One year later, a paper by Zhu \emph{et al.} \cite{zhu_tie_01} studied the observability of non-salient (surface) permanent magnet synchronous machine (SPMSM), using the same aforementioned observability approach. It is shown that the only critical situation for the SPMSM is the standstill operating condition. Other papers that deal with local weak observability conditions of both IM \cite{ghanes_cdc_06} and SPMSM \cite{ezzat_cdc_10,zgorski_sled_12} were published later, resulting in the same observability conditions. Other approaches, such as global observability \cite{rojas_automatica_04}, differential algebraic observability \cite{chiasson_tac_06, diallo_tec_15} and energy-based observability analysis \cite{basic_tac_09}, were proposed, however the results are less promising than those of local observability approach \cite{koteich_tec_15}.  

The observability of the interior (I) permanent magnet synchronous machine (PMSM) is studied in \cite{zaltni_cdc_10, vaclavek_tie_13, hamida_tcon_14}, however, the results of these papers are inconsistent with each other. Furthermore, some of these results are not correctly formulated \cite{koteich_tie_15}. This is due to the fact that the rotor saliency of the IPMSM yields a highly nonlinear model.

For the wound rotor synchronous machine (WRSM) and the synchronous reluctance machine (SyRM), no observability studies were found in the literature \cite{koteich_ecc_15}.

\subsection{Paper contribution}

Based on the above literature review, the following paper presents several contributions. First of all, both SM and IM are studied using the local observability theory. The SMs are studied using a unified approach, where the WRSM is considered as the general case of SMs. This results in a unified observability conditions formulation by introducing the \emph{observability vector} concept for SMs. Thus, the observability conditions of the IPMSM and the SyRM are presented, and the previous results for the SPMSM are found again. Based on WRSM observability analysis, a high-frequency (HF) injection based technique is proposed, combined with observer-based estimation, to guarantee the drive observability for all operating conditions. Numerical simulations and experimental data analysis are carried out to prove the validity of the results. The observability conditions for the IM are formulated in a way to include all previous literature results and are validated via numerical simulations and based on experimental data.

After this introduction, some preliminaries on observer design and observability analysis of dynamical systems are recalled in section 2, and the DC machines observability analysis is presented as an introductory example. The local observability of SMs and IMs is studied in sections 3 and 4 respectively, simulation and experimental results are presented in sections 5 and 6 respectively.

Throughout this paper, AC machines' models are expressed in stationary two-phase reference frame, with conventional modeling assumptions (lossless non-saturated magnetic circuit, sinusoidally distributed magneto-motive forces, constant resistances, etc.). The machines input, output and parameters are considered to be well known. Symbolic math software is used to perform complex calculation.
\section{Observability analysis of dynamical systems}

Nonlinear systems of the following general form (denoted $\Sigma$) are considered:
\begin{subequations}
\begin{eqnarray}
\dot{x} &=& f\left(x, u\right)\\
y &=& h\left(x\right)
\end{eqnarray}
\label{sigma}
\end{subequations}
where $x \in \mathbb{R}^n$ is the state vector, $u \in \mathbb{R}^m$ is the input vector (control signals), $y \in \mathbb{R}^p $ is the output vector (measurement signals), $f$ and $h$ are $C^\infty$ functions. A linear time-invariant (LTI) system is a particular case of $\Sigma$, it has the following form:
\begin{subequations}
\begin{eqnarray}
\dot{x} &=& \mathbf{A} x + \mathbf{B} u \\
y &=& \mathbf{C} x
\end{eqnarray}
\label{sigma_l}
\end{subequations}
where $\mathbf{A}$, $\mathbf{B}$ and $\mathbf{C}$ are state, input and output matrices with appropriate dimensions.

A state observer is a model-based estimator that uses the measurement for correction; it reconstructs the state $x$ of the system from the knowledge of its input $u$ and output $y$. The general form of an observer for the system $\Sigma$ is:
\begin{eqnarray}
\dot{\hat{x}} = f(\hat{x}, u) + \mathbf{K}(.)\left[y - h(\hat{x})\right]
\label{observer_sigma}
\end{eqnarray}
The matrix $\mathbf{K}$, called \emph{observer gain}, can be constant or time-varying. The first term $f(\hat{x},u)$ is the prediction (or estimation) term, and the second one, $\mathbf{K}(.)\left[y - h(\hat{x})\right]$, is the correction (or innovation) term. For LTI systems, the matrix $\mathbf{K}$ is often constant:
\begin{eqnarray}
\dot{\hat{x}} = \mathbf{A} \hat{x} + \mathbf{B} u + \mathbf{K} \left(y - \mathbf{C} \hat{x}\right)
\end{eqnarray}

Observability of the system $\Sigma$ is an important property to be checked prior to observer design. It tells if the state $x$ can be estimated, unambiguously, from the system's input/output signals.

\subsection{Observer design for dynamical systems}
As mentioned earlier, an observer is a replica of the system's dynamical model, with a correction term function of the output estimation error. Observer design for LTI systems is a mature domain. The state estimation error $\tilde{x} = x - \hat{x}$ has the following dynamics:
\begin{eqnarray}
\dot{\tilde{x}} = \left(\mathbf{A} - \mathbf{K} \mathbf{C}\right) \tilde{x}
\end{eqnarray}
If the system is observable, then the eigenvalues of the matrix $\mathbf{A - KC}$ can be arbitrarily placed in the left-half complex plane, by tuning the observer gain $\mathbf{K}$ \cite{luenberger_tme_64}. Therefore, the equilibrium point $\tilde{x} = 0$ is asymptotically stable.

Observer design for nonlinear systems is a more complicated task. Nonetheless, it is a common practice to rely on the linearized model to design \emph{local} observers for nonlinear systems \cite{khalil}: from equations \eqref{sigma} and \eqref{observer_sigma}, we get the following estimation error dynamics:
\begin{eqnarray}
\dot{\tilde{x}} = f(x,u) - f(\hat{x},u) - \mathbf{K}(.)\left[h(x)-h(\hat{x})\right]
\label{sigma_error}
\end{eqnarray}
The observer gain matrix $\mathbf{K}$ should be designed to stabilize the linearized system at $\tilde{x} = 0$. The linearization of the error dynamics \eqref{sigma_error} gives:
\begin{eqnarray}
\dot{\tilde{x}} = \left[\frac{\partial f}{\partial x} (x,u) - \mathbf{K} \frac{\partial h}{\partial x}(x)\right] \tilde{x}
\label{nl_error}
\end{eqnarray}

\begin{figure}[!ht]
	\vspace{-5pt}
	\centering
	\includegraphics[width = \linewidth]{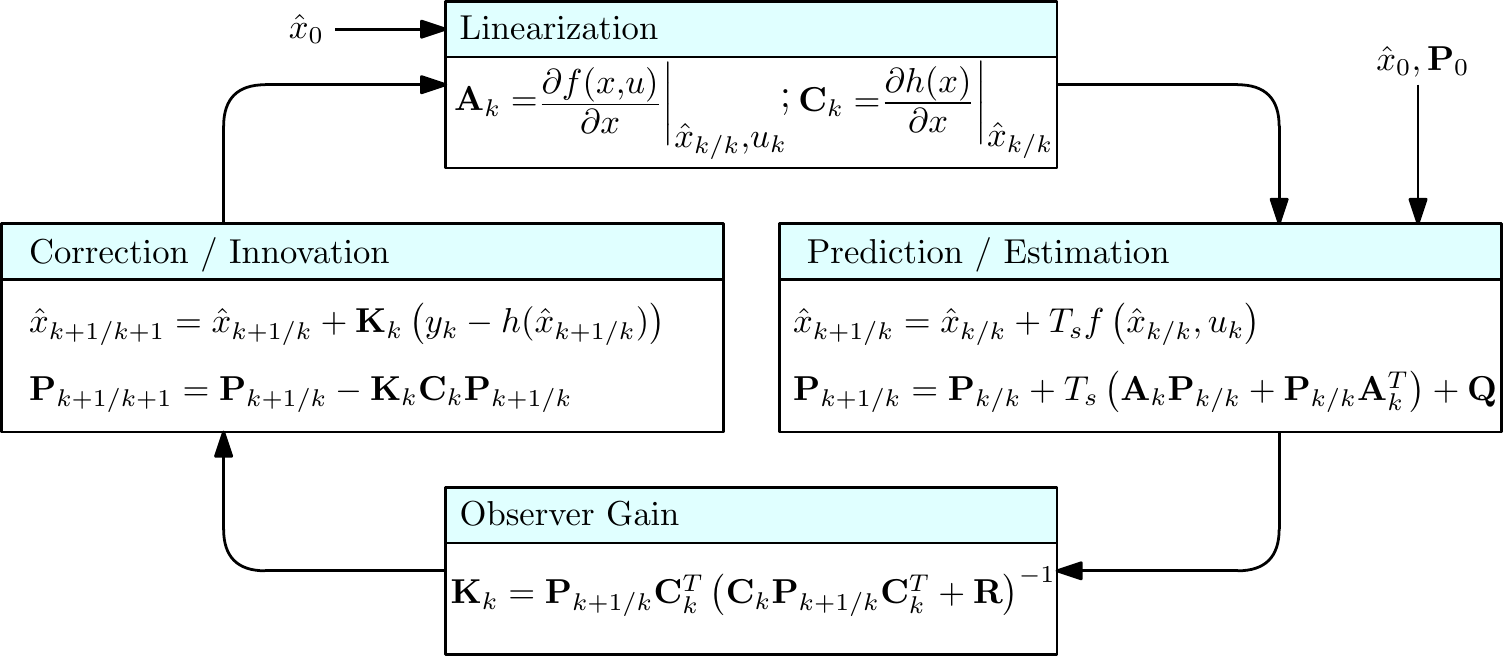}
	\caption{Discrete-time extended Kalman filter}
	\vspace{-15pt}
	\label{ekf}
\end{figure}

The observer gain $\mathbf{K}$ can be designed for a given equilibrium point $x = x_{ss}$ (when $u = u_{ss}$), about which the linearized system is described by:
\begin{eqnarray}
\mathbf{A} = \frac{\partial f}{\partial x} (x_{ss}, u_{ss}) ~~~;~~~~ \mathbf{C} = \frac{\partial h}{\partial x} (x_{ss})
\end{eqnarray}
In this case, $\mathbf{K}$ is constant. However, the observer is guaranteed to work only if $\| \tilde{x}(0) \|$ is sufficiently small and $x$ and $u$ are close enough to their equilibrium values. Yet, the linearization can be taken about the estimate $\hat{x}$, and the observer gain is therefore time varying. A well-know technique to efficiently calculate the matrix $\mathbf{K}$ in real time is the Kalman filter algorithm \cite{kalman_asme_60}, in its deterministic extended version. A sufficiently small initial estimation error $\| \tilde{x}(0) \|$ is still required, but the state $x$ and the input $u$ can be any well-defined trajectories that have no finite escape time \cite{khalil}. The discrete-time extended Kalman filter (EKF) is widely used in industry applications, its algorithm is described in the figure \ref{ekf}. It can be noticed that, for the EKF, only real-time local observability is required.

The matrices $\mathbf{Q}$ and $\mathbf{R}$ in the EKF algorithm are symmetric and positive definite. They are used for tuning of the observer dynamics: higher values of $\mathbf{Q}$ coefficients with respect to $\mathbf{R}$ coefficients result in faster (and more noisy) dynamics, whereas higher values of $\mathbf{R}$ coefficients relatively to $\mathbf{Q}$ coefficients result in smoother, yet slower, observer dynamics. There exist no systematic method to calculate the coefficients of $\mathbf{Q}$ and $\mathbf{R}$, and trial-and-error method is often used. $\mathbf{P}_0$ is a symmetric and positive definite matrix that affects the observer's behavior only at the system start-up. It reflects the expected initial estimation error covariance.

\subsection{Observability analysis}
Let $x_0$ and $x_1$ be two distinct initial states of the system $\Sigma$ at the time $t_0$ ($x_0, x_1 \in \mathbb{R}^n$).
\begin{definition}[Indistinguishability]
The pair $(x_0, x_1)$ is said to be indistinguishable if, for any admissible input $u(t)$, the system outputs $y_0(t)$ and $y_1(t)$, respectively associated to $x_0$ and $x_1$, follow the same trajectories from $t_0$ to $t$ ($\forall t>t_0$), i.e. starting from those two initial states, the system realizes the same input-output map \cite{hermann_tac_77}. 
\end{definition}
In the sequel, $\mathfrak{I}(x_0)$ denotes the set of indistinguishable points from $x_0$.

\begin{definition}[Observability]
	A dynamical system \eqref{sigma} is said to be observable at $x_0$ if $\mathfrak{I}(x_0) = \{x_0\}$. Furthermore, a system $\Sigma$ is observable if it does not admit any indistinguishable pair of states.
\end{definition}

Observability of linear systems can be verified by applying the Kalman criterion: an LTI system is observable if and only if the following observability matrix is full rank:
\begin{eqnarray}
\mathcal{O}_y = \begin{bmatrix}
\mathbf{C}^T & \mathbf{A}^T \mathbf{C}^T & \mathbf{A}^{2^T} \mathbf{C}^T & \ldots & \mathbf{A}^{{n-1}^T} \mathbf{C}^T
\end{bmatrix}^T \nonumber
\end{eqnarray}

For nonlinear systems, the previous definition is rather global. For instance, when applying EKF, only local short-time observability matters. Thus, distinguishing a state from its neighbors is quite enough.

\begin{definition}[$\mathbb{U}$-indistinguishability]
Let $\mathbb{U} \subset \mathbb{R}^n$ be a neighborhood of $x_0$ and $x_1$. $x_0$ and $x_1$ are said to be $\mathbb{U}$-indistinguishable if, for any admissible input $u(t)$, they are indistinguishable when considering time intervals for
which trajectories remain in $\mathbb{U}$. The set of states that are $\mathbb{U}$-indistinguishable of $x_0$ is denoted $\mathfrak{I}_{\mathbb{U}}(x_0)$.
\end{definition}

\begin{definition}[Local observability]
A system is locally observable (LO) at $x_0$ if, for any neighborhood $\mathbb{U}$ of $x_0$, $\mathfrak{I}_{\mathbb{U}}(x_0) = \{x_0\}$. If this is true for any $x_0$ the system is then locally observable.
\end{definition}

\begin{definition}[Local weak observability]
A system \eqref{sigma} is said to be locally weakly observable (LWO) if, for any $x_0$, there exists a neighborhood $\mathbb{V}$ of $x_0$ such that for any neighborhood $\mathbb{U}$ of $x_0$ contained in $\mathbb{V}$ ($\mathbb{U} \subset \mathbb{V}$) , there is no $\mathbb{U}$-indistinguishable state from $x_0$. Equivalently, a system is LWO if:
\begin{equation}
\forall x_0, ~\exists \mathbb{V}, ~ x_0 \in \mathbb{V}, ~ \forall \mathbb{U} \subset \mathbb{V}, ~ \mathfrak{I}_{\mathbb{U}}(x_0)=\{x_0\} \nonumber
\end{equation}
\end{definition}

Roughly speaking, if one can instantaneously distinguish any state of a system from its neighbors, the system is LWO.

\begin{definition}[Observability rank condition]
The system $\Sigma$ is said to satisfy the observability rank condition at $x_0$, if the observability matrix, denoted by $\mathcal{O}_y(x)$, is full rank at $x_0$. $\mathcal{O}_y(x)$ is given by:
\begin{equation}
\mathcal{O}_y(x) = \frac{\partial}{\partial x}\left[ \begin{matrix}
\mathcal{L}^0_fh(x)^T & 
\mathcal{L}_fh(x)^T & 
\mathcal{L}_f^2h(x)^T & 
\ldots & 
\mathcal{L}_f^{n-1}h(x)^T
\end{matrix} \right]_{x=x_0}^T \nonumber
\end{equation}
where $\mathcal{L}_f^{k}h(x)$ is the $k$th-order \emph{Lie derivative} of the function $h$ with respect to the vector field $f$.
\end{definition}

%

\begin{theorem}[\cite{hermann_tac_77}]A system $\Sigma$ satisfying the observability rank condition at $x_0$ is locally weakly observable at $x_0$. More generally, a system $\Sigma$ satisfying the observability rank condition, for any $x_0$, is LWO.
\end{theorem}

\begin{remark}
	In the case of LTI systems, rank criterion is equivalent to Kalman criterion. However, locality is meaningless for such systems.
\end{remark}

\begin{remark}
Rank criterion gives only a sufficient observability condition for nonlinear systems.
\end{remark}

\subsection{Introductory example: DC machines}
In order to introduce the reader to the observability theory application for sensorless electric drives, two DC machines are studied at first: permanent magnet (PM-) DC machines, which has a linear model, and series (S-) DCM which has a nonlinear model.

For both machines, the observability of the rotor speed $\Omega$ and the identifiability\footnote{We call identifiability the observability of (supposed) constant states.} of the load torque $T_l$ will be studied. The only measurement is the armature current $i_a$ and the input voltage $v$ is considered to be known.

\subsubsection{PM-DCM observability}
The state-space model of a PM-DCM can be written as follows:
\begin{subequations}
\begin{eqnarray*}
\frac{di_a}{dt} &=& \frac{1}{L_a} \left( v - R_a i_a - K_e~\Omega \right)\\
\frac{d \Omega}{dt} &=& \frac{1}{J} \left(K_e~i_a - T_l\right)- \frac{f_v}{J} \Omega\\
\frac{d T_l}{dt} &=& 0
\end{eqnarray*} 
\end{subequations}
It is a linear model, the application of Kalman criterion gives the following observability matrix:
\begin{eqnarray}
\mathcal{O}_y^{PMDC} = \begin{bmatrix}
\mathbf{C} \\ \mathbf{CA} \\ \mathbf{CA}^2
\end{bmatrix} = \begin{bmatrix}
1 & 0 & 0\\
-\frac{R_a}{L_a} & -\frac{K_e}{L_a} & 0\\
\frac{R_a^2}{L_a^2} - \frac{K_e^2}{J L_a} & \frac{K_e}{L_a} \left(\frac{R_a}{L_a} + \frac{f_v}{J}\right) & \frac{K_e}{J L_a}
\end{bmatrix} \nonumber
\end{eqnarray}
The determinant of this matrix is:
\begin{eqnarray}
\Delta_y^{PMDC}  = -\frac{K_e^2}{J L_a^2} \neq 0 \nonumber
\end{eqnarray}
Thus, the PM-DCM is observable as far as $K_e$, proportional to the field PM flux, is non-zero.

\subsubsection{S-DCM observability}
In an S-DCM, the field coil is connected in series with the armature. This results in the following S-DCM nonlinear state-space model: 
\begin{subequations}
\begin{eqnarray*}
\frac{d i_a}{dt} &=& \frac{1}{L_a + L_f} \left( v - (R_a + R_f) i_a - K i_a \Omega \right)\\
\frac{d \Omega}{dt} &=& \frac{1}{J} \left(K i_a^2 - T_l\right) - \frac{f_v}{J} \Omega\\
\frac{d T_l}{dt} &=& 0
\end{eqnarray*}
\end{subequations}
The output and its derivatives are:
\begin{eqnarray*}
h(x) &=& i_a \\
\mathcal{L}_fh(x) &=& \frac{1}{L} \left(v - R i_a - K \Omega i_a \right)\\
\mathcal{L}^2_fh(x) &=& -\frac{R + K \Omega}{L^2}\left(v - R i_a - K \Omega i_a \right) - \frac{K i_a}{J L} \left(K i_a^2 - T_l - f_v \Omega \right)
\end{eqnarray*}
where $L = L_a + L_f$ and $R = R_a + R_f$. The observability matrix can be the written as:
\begin{eqnarray*}
\mathcal{O}_y^{SDC}  = \begin{bmatrix}
1 & 0 & 0\\
-\frac{(R + K \Omega)}{L} & -\frac{K i_a}{L} & 0\\
\frac{(R + K \Omega)^2}{L^2} - \frac{K (3 K i_a^2 - T_l)}{J L} & -\frac{K (v - 2 i_a (R + K \Omega))}{L^2}  + \frac{K f_v i_a}{J L} & \frac{K i_a}{J L} 
\end{bmatrix}
\end{eqnarray*}
It has the following determinant:
\begin{eqnarray*}
\Delta_y^{SDC} = -\frac{K^2}{J L^2} i_a^2
\end{eqnarray*}
Therefore, the observability of an S-DCM cannot be guaranteed if the current $i_a$ is zero. Indeed, if the current is null, the back-EMF is null too, and no information about the speed and the load torque can be found in the voltage equation:
\begin{eqnarray*}
\frac{di_a}{dt} = \frac{1}{L} v
\end{eqnarray*}
The previous results can be generalized for other DCMs: the DCM observability is not guaranteed in absence of field magnetic flux.
\section{Observability analysis of synchronous machines}
A synchronous machine (SM) is an AC machine. In motor operating mode, the stator is fed by a polyphase voltage to generate a rotating magnetic field. Depending on the rotor structure, the rotating field interacts with the rotor in different ways. There exist, broadly, three types of interaction between the stator magnetic field and the rotor: 1) magnetic field interaction with electromagnet, which is the case of the WRSM, 2) magnetic field with permanent magnet, in the case of the PMSM and 3) magnetic field with ferromagnetic material, which is the case of SyRM, where the salient type rotor moves in a way to provide a minimum reluctance magnetic path for the stator rotating field.

These types of interactions can be combined. For instance, a SyRM can be assisted with PMs in order to increase its electromagnetic torque. Moreover, the rotor of a WRSM or a PMSM can be designed in a way to create an anisotropic magnetic path (salient type rotors) to improve the machine performances. In this paper, interior (I) PMSM stands for any brushless SM with a salient rotor having PMs. The salient type WRSM is considered to be the most general SM, and it will be shown that the extension of the study of this machine to the other SMs is not a real burden.

In sensorless synchronous drives, currents are measured, voltages are considered to be known, and rotor position and angular speed are to be estimated. 
In this section, the observability of salient type WRSM is studied, which results in the definition of a new concept, the \emph{observability vector}, that unifies the observability conditions formulation for SMs.

\begin{figure}[!ht]
	\centering
\includegraphics[scale = 0.85]{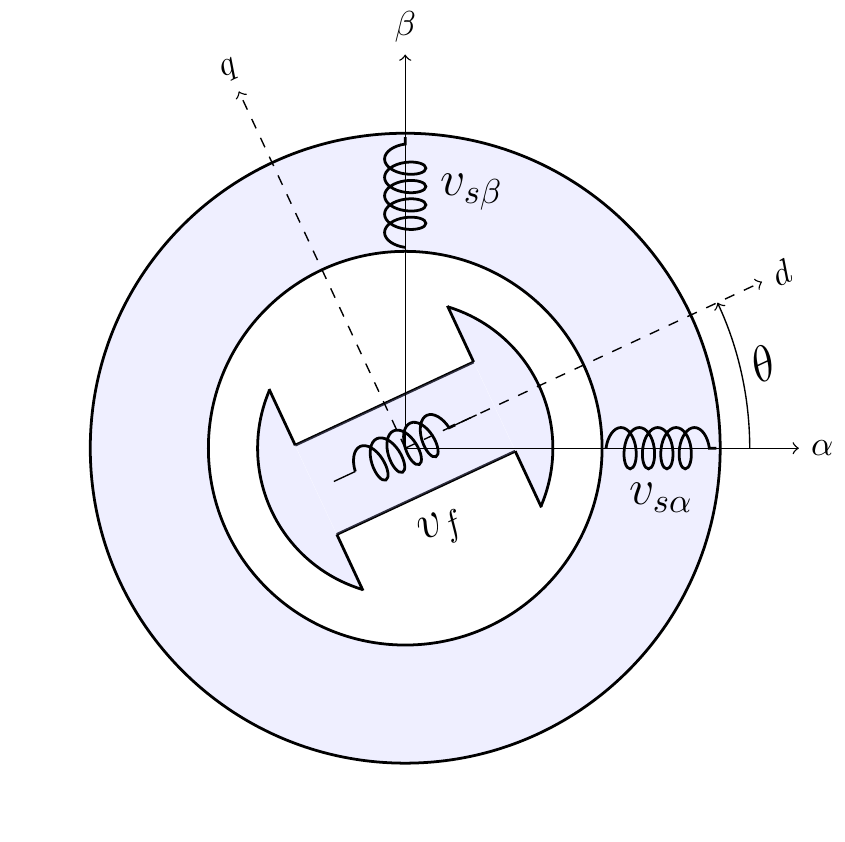}
\caption{Schematic representation of the WRSM in the $\alpha \beta$ reference frame }
\label{wrsm}
\end{figure}

\subsection{WRSM model} 
In the view of observability analysis, it is preferred to express state variables in the reference frame where measurement is performed. For this reason, the state-space model of the WRSM is written in the two-phase ($\alpha \beta$) stationary reference frame for the stator currents (see Figure \ref{wrsm}), and in the rotor reference frame ($dq$) for the rotor current ($i_f$). The speed dynamics are neglected for two reasons: 1) they are slower than the currents dynamics and 2) there is no slip phenomenon in the SM. This gives the following model:
\begin{subequations}
\begin{eqnarray}
\frac{d{\mathcal{I}}}{dt}&=&-{\mathfrak{L}(\theta)^{-1}}{\mathfrak{R}_{eq}}\mathcal{I}
+ {\mathfrak{L}(\theta)^{-1}}\mathcal{V}\\
\frac{d\omega}{dt} &=& 0 \\ 
\frac{d\theta}{dt} &=& \omega
\end{eqnarray}
\label{ss_sm}
\end{subequations}
This model can be fitted to the structure \eqref{sigma}, where the state, input and output vectors are respectively:
\begin{eqnarray*}
x = \left[
\begin{matrix}
\mathcal{I}^T & \omega & \theta
\end{matrix}
\right]^T~~~~~;~~~~~
u = 
\mathcal{V}= \left[\begin{matrix}
v_{s \alpha} & v_{s \beta} & v_f
\end{matrix}\right]^T~~~~~;~~~~~
y = 
\mathcal{I} = \left[\begin{matrix}
i_{s \alpha} & 
i_{s \beta} & 
i_f
\end{matrix}\right]^T
\end{eqnarray*}
$\mathcal{I}$ and $\mathcal{V}$ are the machine current and voltage vectors. Indices $\alpha$ and $\beta$ stand for stator signals, index $f$ stands for rotor (field) ones. $\omega$ stands for the electrical rotor angular speed and $\theta$ for the electrical angular position of the rotor. $\mathfrak{L}$ is the (position-dependent) matrix of inductances, and $\mathfrak{R}_{eq}$ is the equivalent resistance matrix that extends the resistance matrix $\mathfrak{R}$:
\begin{equation*}
\mathfrak{L}=\left[ \begin{matrix}
   {{L}_{0}}+{{L}_{2}}\cos 2\theta  & {{L}_{2}}\sin 2\theta  & {{M}_{f}}\cos \theta   \\
   {{L}_{2}}\sin 2\theta  & {{L}_{0}}-{{L}_{2}}\cos 2\theta  & {{M}_{f}}\sin \theta   \\
   {{M}_{f}}\cos \theta  & {{M}_{f}}\sin \theta  & {{L}_{f}}  \\
\end{matrix} \right] ~;~
\mathfrak{R} = \left[\begin{matrix}
R_s & 0   &  0 \\
  0 & R_s &  0 \\
  0 & 0   & R_f
\end{matrix}\right] ~;~
\mathfrak{R}_{eq} = \mathfrak{R} + \frac{\partial\mathfrak{L}}{\partial\theta}\omega
\end{equation*}

\subsection{WRSM observability}  The system \eqref{ss_sm} is a 5-${th}$ order system. Its observability study requires the evaluation of the output derivatives up to the 4-${th}$ order. In this study, only the first order derivatives are calculated, higher order ones are very difficult to calculate and to deal with. This gives the following ``partial'' observability matrix:
\begin{equation}
\mathcal{O}_{y}^{SM} = \left[ \begin{matrix}
\mathbf{I}_{3} &\mathbf{O}_{3 \times 1} &\mathbf{O}_{3 \times 1}  \\
-\mathfrak{L}^{-1}\mathfrak{R}_{eq} &-\mathfrak{L}^{-1}\mathfrak{L}' \mathcal{I} &\mathfrak{L^{-1}}' \mathfrak{L} \frac{d \mathcal{I}}{dt} - \mathfrak{L^{-1}}\mathfrak{L}''\omega \mathcal{I} 
\end{matrix} \right] 
\label{obsv_matrix}
\end{equation}
where $\mathbf{I}_n$ is the $n \times n$ identity matrix and $\mathbf{O}_{n \times m}$ is an $n \times m$ zero matrix.
$\mathfrak{L}'$ and $\mathfrak{L}''$ denote, respectively, the first and second partial derivatives of $\mathfrak{L}$ with respect to $\theta$.


It is sufficient to have five linearly independent lines of matrix \eqref{obsv_matrix} to ensure the local observability of the system.
The first five lines that come from the currents and the first derivatives of $i_{s \alpha}$ and $i_{s \beta}$ are studied. This choice is motivated by the fact that the stator currents are available for measurement in all synchronous machines, the rotor current (that gives the sixth line of matrix \eqref{obsv_matrix}) does not exist in the PMSM and the SyRM. Another reason comes from the physics of the machine: $i_f$ is a DC signal, whereas both $i_{s \alpha}$ and $i_{s \beta}$ are AC signals, so it is more convenient for physical interpretation to take them together.

Currents are expressed in the rotor ($dq$) reference frame, using the following Park transformation, in order to make the interpretation easier:
\begin{eqnarray*}
\begin{bmatrix}
i_{s \alpha} \\ i_{s \beta}
\end{bmatrix}
 =
\begin{bmatrix}
\cos \theta & -\sin \theta \\
\sin \theta & \cos \theta
\end{bmatrix}
\begin{bmatrix}
i_{sd} \\ i_{sq}
\end{bmatrix}
\end{eqnarray*}
The determinant is written, in function of stator currents and fluxes, under the following general form:
\begin{eqnarray}
\Delta_y^{WRSM} =  \mathcal{D}^{WRSM}\omega + \mathcal{N}^{WRSM}
\label{wrsm_determinant}
\end{eqnarray}
where
\begin{eqnarray}
\mathcal{D}^{WRSM} &=&  \frac{1}{\sigma_d L_d L_q} \left[
\left(\psi_{sd} - L_q i_{sd}\right)^2 + \sigma_\Delta L_\Delta^2 i_q^2
\right] \label{wrsm_determinant_d}\\
\mathcal{N}^{WRSM} &=& \frac{\sigma_\Delta}{\sigma_d}\frac{ L_\Delta}{L_d L_q} \left[
\frac{d\psi_{sd}}{dt}i_{sq} + \frac{d\psi_{sq}}{dt}i_{sd} - \left(\frac{di_{sd}}{dt} \psi_{sq} + \frac{di_{sq}}{dt} \psi_{sd}  \right) \right] \label{wrsm_determinant_n}
\end{eqnarray}
with
\begin{eqnarray*}
L_d = L_0 + L_2~~;~~ L_q = L_0 - L_2~~;~~ L_\Delta = L_d - L_q = 2 L_2
\end{eqnarray*}
and:
\begin{eqnarray}
\sigma_d = 1 - \frac{M_f^2}{L_d L_f} ~~~~;~~~~ \sigma_\Delta = 1 -  \frac{M_f^2}{L_\Delta L_f}
\end{eqnarray}
Note that the stator fluxes can be expressed as follows:
\begin{eqnarray}
\psi_{sd} &=& L_d i_{sd} + M_f i_f \\
\psi_{sq} &=& L_q i_{sq}
\end{eqnarray}
It is sufficient for the determinant \eqref{wrsm_determinant} to be non-zero to guarantee the machine's observability. This condition is examined below for each type of synchronous machine.

\subsection{Brushless SMs models and observability conditions}
\begin{figure}[!ht]
\centering
\begin{subfloat}[IPMSM]{
  \centering
  \includegraphics[scale = 0.6]{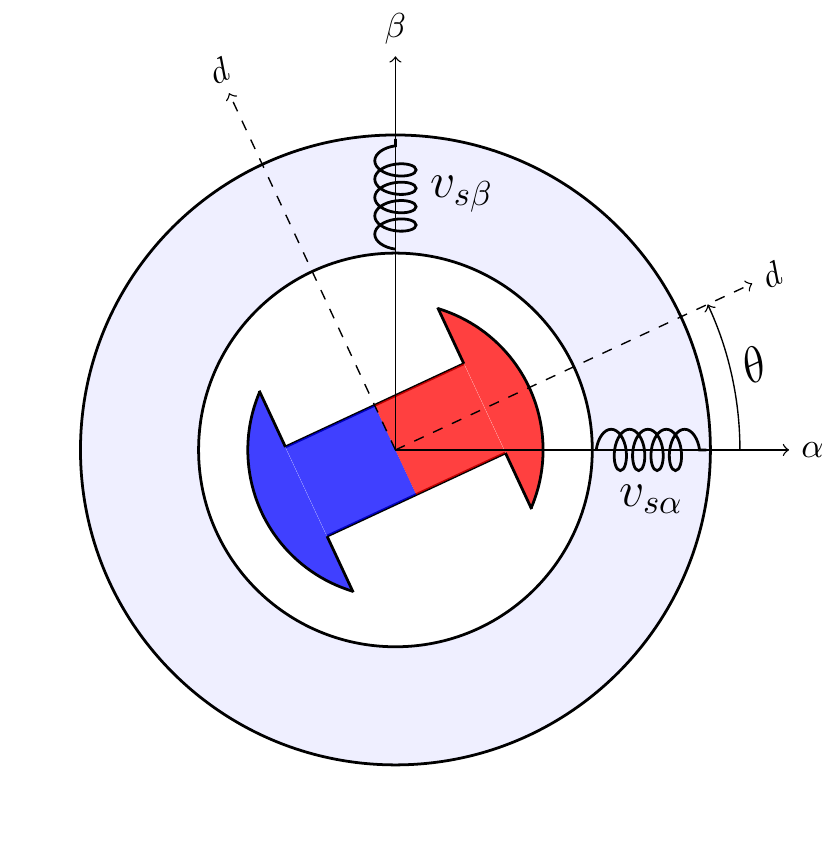}
  \label{ipmsm_schema}}
\end{subfloat}%
\begin{subfloat}[SPMSM]{
  \centering
  \includegraphics[scale = 0.6]{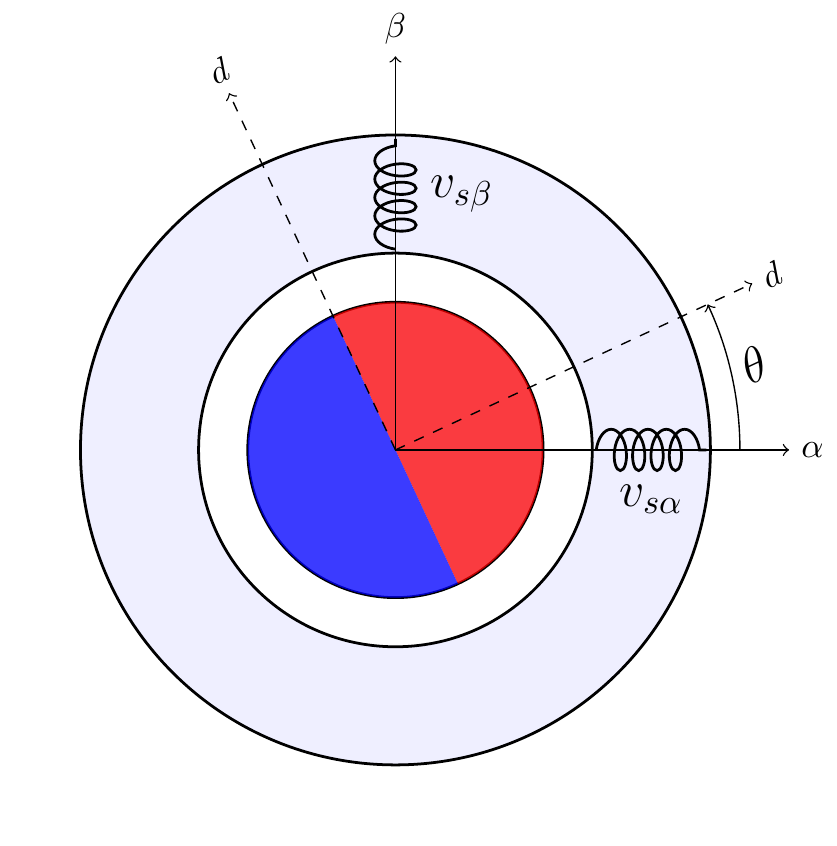}
  \label{spmsm_schema}}
\end{subfloat}
\begin{subfloat}[SyRM]{
\includegraphics[scale = 0.6]{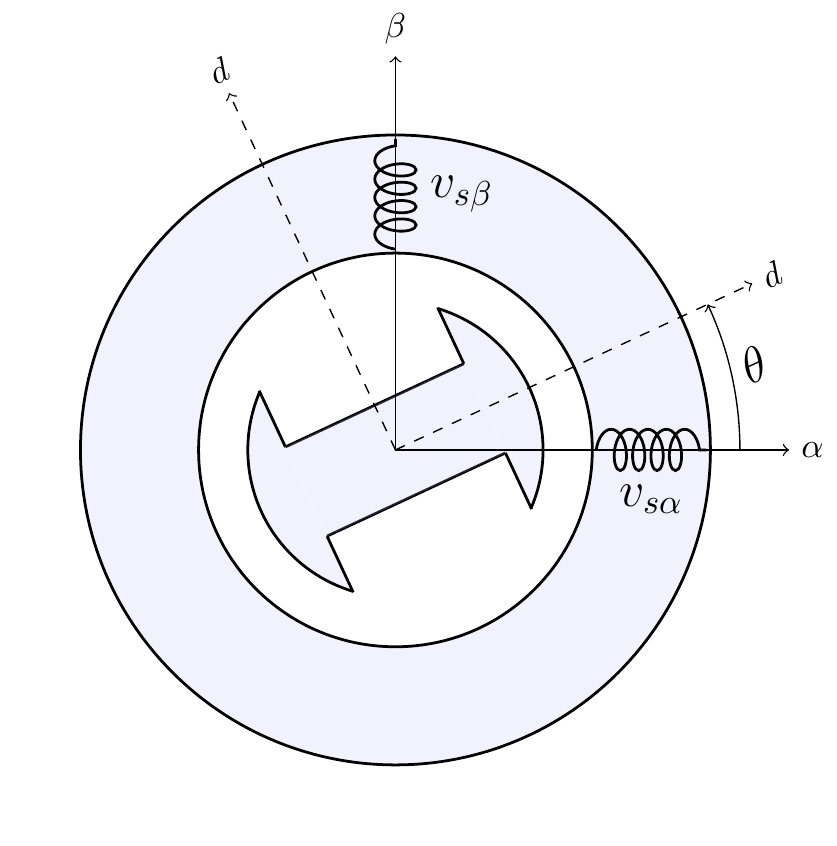}
\label{syrm_schema}}
\end{subfloat}
\caption{Schematic representation of (a) the IPMSM, (b) the SPMSM and (c) the SyRM}
\label{pmsm}
\end{figure}

The brushless SMs can be seen as special cases of the salient-type WRSM. The IPMSM (Figure \ref{ipmsm_schema}) model can be derived from the WRSM one by considering the rotor magnetic flux to be constant:
\begin{eqnarray}
\frac{d\psi_f}{dt} = 0
\label{cond_pmsm_1}
\end{eqnarray}
and by substituting $M_f i_f$ of the WRSM by the permanent magnet flux $\psi_r$ of the IPMSM:
\begin{eqnarray}
i_f = \frac{\psi_r}{M_f}
\label{cond_pmsm_2}
\end{eqnarray}

The equations of SPMSM (Figure \ref{spmsm_schema}) are the same as the IPMSM, except that the stator self-inductances are constant and independent of the rotor position, that is:
\begin{eqnarray}
L_2 = 0 ~~~\implies~~~ L_d = L_q = L_0
\label{non_salient}
\end{eqnarray}

The SyRM (Figure \ref{syrm_schema}) model can be derived from the IPMSM model by considering the rotor magnetic flux $\psi_r$ to be zero: 
\begin{eqnarray}
\psi_r \equiv 0
\label{cond_syrm}
\end{eqnarray}

Consequently, the determinant of the observability matrix of brushless synchronous machines can be derived from the equations \eqref{wrsm_determinant}, \eqref{wrsm_determinant_n} and \eqref{wrsm_determinant_d} by substituting:
\begin{eqnarray*}
\sigma_d = \sigma_\delta = 1
\end{eqnarray*}
and for the SPMSM by substituting in addition:
\begin{eqnarray*}
L_\Delta = L_d - L_q = 0
\end{eqnarray*}

This gives the following determinant for the IPMSM and SyRM, $$\Delta_y^{IPMSM} = \mathcal{D}^{IPMSM} \omega + \mathcal{N}^{IPMSM}$$ where:
\begin{eqnarray}
\mathcal{D}^{IPMSM} &=&  \frac{1}{L_d L_q} \left[
\left(\psi_{sd} - L_q i_{sd}\right)^2 + L_\Delta^2 i_q^2
\right] \label{ipmsm_determinant_d}\\
\mathcal{N}^{IPMSM} &=& \frac{ L_\Delta}{L_d L_q} \left[
\frac{d\psi_{sd}}{dt}i_{sq} + \frac{d\psi_{sq}}{dt}i_{sd} - \left(\frac{di_{sd}}{dt} \psi_{sq} + \frac{di_{sq}}{dt} \psi_{sd}  \right) \right] \label{ipmsm_determinant_n}
\end{eqnarray}
and for the SPMSM :
\begin{eqnarray}
\Delta_y^{SPMSM} &=&  \frac{1}{L_0^2} \left[
\left(\psi_{sd} - L_0 i_{sd}\right)^2
\right] \omega \label{spmsm_determinant}
\end{eqnarray}
The observability of each machine is guaranteed as far as the determinant of its observability matrix is different from zero.

\subsection{The concept of observability vector}
In this paragraph, the concept of \emph{observability vector} is introduced for SMs. This concept allows us to formulate a unified observability condition, that can be tested in real time, for all sensorless synchronous drives. The observability vector, $\Psi_\mathcal{O}$, is defined as follows for the WRSM:
\begin{eqnarray*}
\Psi_{\mathcal{O} d} &=& L_\Delta i_{sd} + M_f i_f\\
\Psi_{\mathcal{O} q} &=& \sigma_\Delta L_\Delta i_{sq}
\end{eqnarray*}
Therefore, it can be expressed for the other SMs:
\begin{eqnarray*}
\text{IPMSM} &:& \Psi_{\mathcal{O} d} = L_\Delta i_{sd} + \psi_r ~~;~~ \Psi_{\mathcal{O} q} = L_\Delta i_{sq} \\
\text{SPMSM} &:& \Psi_{\mathcal{O} d} = \psi_r ~~~~~~~~~~~~~;~~  \Psi_{\mathcal{O} q} = 0 \\
\text{SyRM} &:& \Psi_{\mathcal{O} d} = L_\Delta i_{sd} ~~~~~~~~~;~~  \Psi_{\mathcal{O} q} = L_\Delta i_{sq}
\end{eqnarray*}

\begin{figure}[!ht]
	\centering
	\includegraphics[scale=1.2]{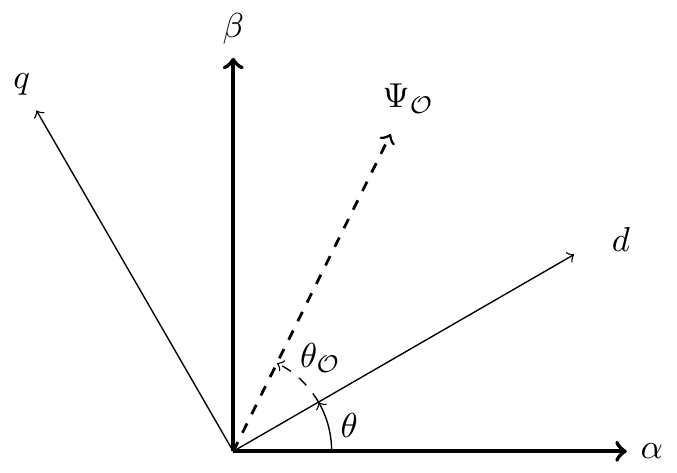}
	\caption{Observability vector}%
	\label{obsv_vect}
\end{figure}

The unified observability condition can be formulated as follows (figure \ref{obsv_vect}):
\begin{proposition}
A synchronous machine is locally observable if the angular velocity of the observability vector in the rotor reference frame is different from the electrical velocity of the rotor in the stationary reference frame:
\begin{eqnarray*}
\frac{d \theta_\mathcal{O}}{dt} \neq \omega
\end{eqnarray*}
\end{proposition}

\begin{proof}
After substituting the stator fluxes by their equations $\psi_{sd} = L_d i_{sd} + M_f i_f$ and $\psi_{sq} = L_q i_{sq}$, the determinant of the WRSM observability matrix can be expressed as follows:
\begin{eqnarray}
\Delta^{WRSM} &=&   \frac{1}{\sigma_d {L}_d L_q} \left[
\left(L_\Delta i_{sd} + M_f i_f \right)^2 + \sigma_\Delta L_\Delta^2 i_{sq}^2
\right]\omega \nonumber\\
&& + \frac{\sigma_\Delta}{\sigma_d}\frac{L_\Delta}{L_d L_q} \left[
\left(L_\Delta \frac{di_{sd}}{dt} + M_f \frac{di_f}{dt} \right) i_{sq}  -  \left(L_\Delta i_{sd} + M_f i_f \right) \frac{di_{sq}}{dt} 
\right] 
\label{delta_wrsm} 
\end{eqnarray}

The observability condition $\Delta_{WRSM} \neq 0$ is equivalent to:
\begin{eqnarray*}
\omega \neq  
	\frac{
	\left(L_\Delta i_{sd} + M_f i_f \right) \sigma_\Delta L_\Delta \frac{di_{sq}}{dt} - \left(L_\Delta \frac{di_{sd}}{dt} + M_f \frac{di_f}{dt} \right) \sigma_\Delta L_\Delta i_{sq}
	}
	{
	\left(L_\Delta i_{sd} + M_f i_f \right)^2 + \sigma_\Delta L_\Delta^2 i_{sq}^2
	}
\end{eqnarray*}
The above equation can be rearranged to get the following one:
\begin{eqnarray*}
\omega &\neq& 
	\frac{
		  \left(L_\Delta i_{sd} + M_f i_f \right)^2 + \sigma_\Delta^2 L_\Delta^2 i_{sq}^2
		}
	{
		\left(L_\Delta i_{sd} + M_f i_f \right)^2 + \sigma_\Delta L_\Delta^2 i_{sq}^2
	} \times \nonumber \\
	&& 
	\frac{
		\left(L_\Delta i_{sd} + M_f i_f \right) \sigma_\Delta L_\Delta \frac{di_{sq}}{dt} - \left(L_\Delta \frac{di_{sd}}{dt} + M_f \frac{di_f}{dt} \right) \sigma_\Delta L_\Delta i_{sq}
	}
	{
		\left(L_\Delta i_{sd} + M_f i_f \right)^2 + \sigma_\Delta^2 L_\Delta^2 i_{sq}^2
		} 
\end{eqnarray*}
Then
\begin{eqnarray}
\omega &\neq& 
\frac{
	\left(L_\Delta i_{sd} + M_f i_f \right)^2 + \sigma_\Delta^2 L_\Delta^2 i_{sq}^2
}
{
	\left(L_\Delta i_{sd} + M_f i_f \right)^2 + \sigma_\Delta L_\Delta^2 i_{sq}^2
} \frac{d}{dt} \left(\arctan\frac{\sigma_\Delta L_\Delta i_{sq}}{L_\Delta i_{sd} + M_f i_f}\right) \label{cond_wrsm}
\end{eqnarray}

Let $\Psi_\mathcal{O}$ be a vector having the following components in the $dq$ reference frame:
\begin{eqnarray*}
\Psi_{\mathcal{O} d} &=& L_\Delta i_{sd} + M_f i_f\\
\Psi_{\mathcal{O} q} &=& \sigma_\Delta L_\Delta i_{sq}
\end{eqnarray*}

Let $\theta_\mathcal{O}$ be the angle of this vector in the rotor reference frame, the observability condition can be then written:
\begin{eqnarray*}
\omega \neq 
\left( \frac{
	\left(L_\Delta i_{sd} + M_f i_f \right)^2 + \sigma_\Delta^2 L_\Delta^2 i_{sq}^2
}
{
	\left(L_\Delta i_{sd} + M_f i_f \right)^2 + \sigma_\Delta L_\Delta^2 i_{sq}^2
}
\right) \frac{d \theta_{\mathcal{O}}}{dt}
\end{eqnarray*}

The following approximation is jugged to be reasonable\footnote{Note that this approximation is valid for the WRSM at standstill where $\omega = 0$. Moreover, it is reasonable for non zero speeds. For brushless SMs, where $\sigma_\Delta = 1$, this approximation is an equality.}:
\begin{eqnarray}
\frac{
	\left(L_\Delta i_{sd} + M_f i_f \right)^2 + \sigma_\Delta^2 L_\Delta^2 i_{sq}^2
}
{
	\left(L_\Delta i_{sd} + M_f i_f \right)^2 + \sigma_\Delta L_\Delta^2 i_{sq}^2
} \approx 1
\label{approx_wrsm}
\end{eqnarray}
Therefore, the WRSM is observable if the following equation holds:
\begin{eqnarray}
\omega \neq \omega_{\mathcal{O}}
\end{eqnarray}
where:
\begin{eqnarray}
\omega_{\mathcal{O}} = \frac{d}{dt} \theta_{\mathcal{O}} = \frac{d}{dt} \arctan \left(\frac{\sigma_\Delta L_\Delta i_{sq}}{L_\Delta i_{sd} + M_f i_f}\right)
\end{eqnarray}

The same reasoning applies, without approximation, to the other SMs where:


\begin{eqnarray}
\Delta^{IPMSM} = \frac{1}{L_d L_q} \left[\left(L_\Delta i_{sd} + \psi_r \right)^2 + L_\Delta^2 i_{sq}^2 \right]\omega + \frac{L_\Delta}{L_d L_q} \left[L_\Delta i_{sq} \frac{di_{sd}}{dt} - \left(L_\Delta i_{sd} + \psi_r \right) \frac{di_{sq}}{dt}
\right]
\label{delta_ipmsm}
\end{eqnarray}
\begin{eqnarray}
\Delta^{SPMSM} = \frac{\psi_r^2}{L_0^2}\omega
\label{delta_spmsm}
\end{eqnarray}
\begin{eqnarray}
\Delta^{SyRM} = \frac{L_\Delta^2}{L_d L_q} \left[\left( i_{sd}^2 + i_{sq}^2 \right)\omega + \frac{di_{sd}}{dt} i_{sq} - i_{sd} \frac{di_{sq}}{dt}
\right]
\label{delta_syrm}
\end{eqnarray}

\end{proof}

\begin{figure}[!htb]
	\centering
	\subfloat[WRSM]{\includegraphics[width = 0.45\linewidth]{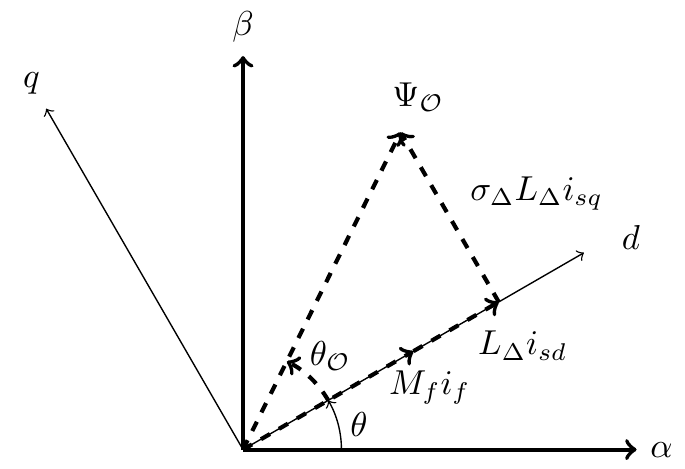}}
	\qquad
	\subfloat[IPMSM]{\includegraphics[width = 0.45\linewidth]{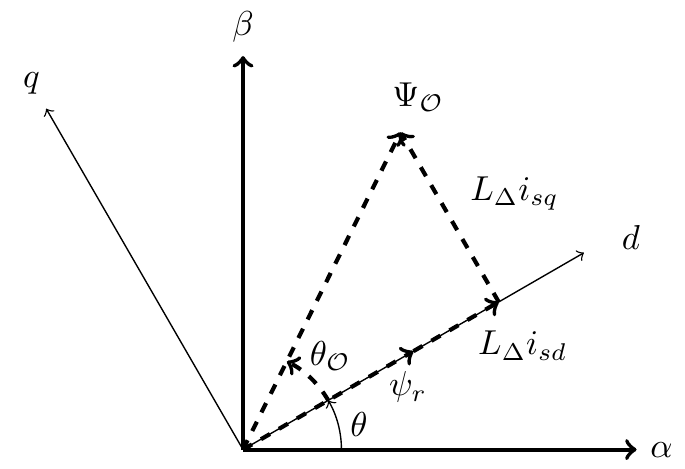}}
	\qquad
	\subfloat[SPMSM]{\includegraphics[width = 0.45\linewidth]{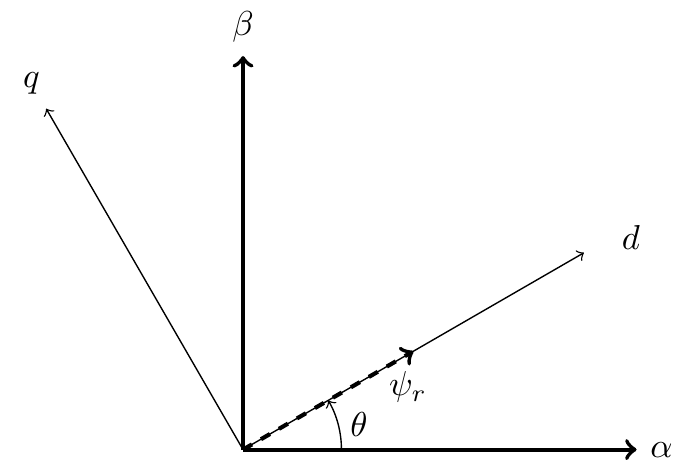}}
	\qquad
	\subfloat[SyRM]{\includegraphics[width = 0.45\linewidth]{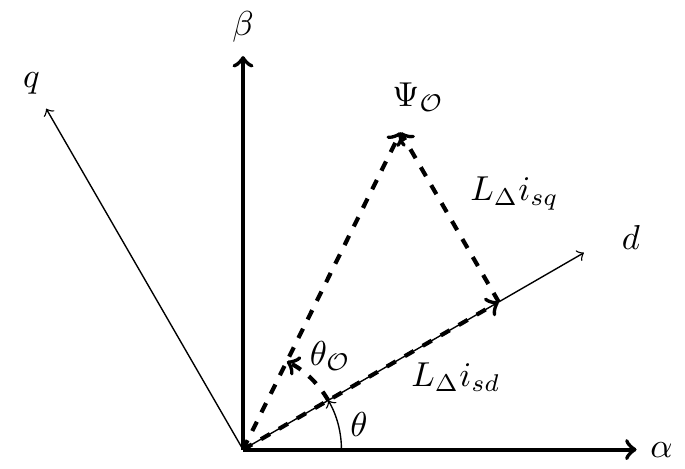}}
	\caption{Observability vector of some synchronous machines}
	\label{obsv_vect_exemples}
\end{figure}

The figure \ref{obsv_vect_exemples} shows the observability vectors of these SMs. Further remarks on the PMSM observability analysis can be found in \cite{koteich_tec_15, koteich_sled_15}.

\begin{remark}
The expression of $\Delta^{WRSM}$ \eqref{delta_wrsm} shows that WRSM can be observable at standstill if the currents $i_{sd}$, $i_{sq}$ and $i_f$ are not constant at the same time. IPMSM \eqref{delta_ipmsm} and SyRM \eqref{delta_syrm} can be also observable at standstill. However, the observability of SPMSM at standstill cannot be guaranteed.
\end{remark}

\begin{remark}
The $d$-axis component of the observability vector is nothing but the so-called \emph{active flux} introduced by \cite{boldea_tpe_08}, which is, by definition, the generalization for all SMs of the flux that multiplies the current $i_{sq}$ in the torque equation of an SPMSM: $T_m = p \psi_r i_{sq}$.
\end{remark}

\begin{remark}
For hybrid excited synchronous machines (HESM), where the wound rotor is assisted with permanent magnets, it can be proved that the $dq$ components of the observability vector are $\Psi_{\mathcal{O}d} = L_\Delta i_{sd} + M_f i_f + \psi_r$ and $\Psi_{\mathcal{O}q} = \sigma_\Delta L_\Delta i_{sq}$.
\end{remark}

The observability conditions formulation with the observability vector concept generalizes, and corrects, some results on the observability of PMSMs \cite{vaclavek_tie_13, koteich_tie_15}. 
\section{Observability analysis of induction machines}
An induction machine is an AC machine. In motor operating mode, the stator is fed by a polyphase voltage to generate a rotating magnetic field, and the rotor is short-circuited. The rotating magnetic field induces currents in the rotor closed circuit, which, in turn, produces a rotor magnetic field. The rotor magnetic field interacts with the stator one to produce a torque that tends to rotate the motor shaft.

This section deals with the local observability conditions of induction machines. The stator currents are measured and the stator voltages are considered to be known. The observability of the IM is studied with and without speed measurement. In the sequel, the speed is not considered to be constant, therefore the resistant torque $T_r$, which includes the load torque and the frictions, is added to the estimated state vector, and is considered to be slowly varying. This resistant torque has a critical impact on the rotor speed, due to slip phenomenon.

\subsection{IM model} 
\begin{figure}[!ht]
	\centering
\includegraphics[scale = 0.85]{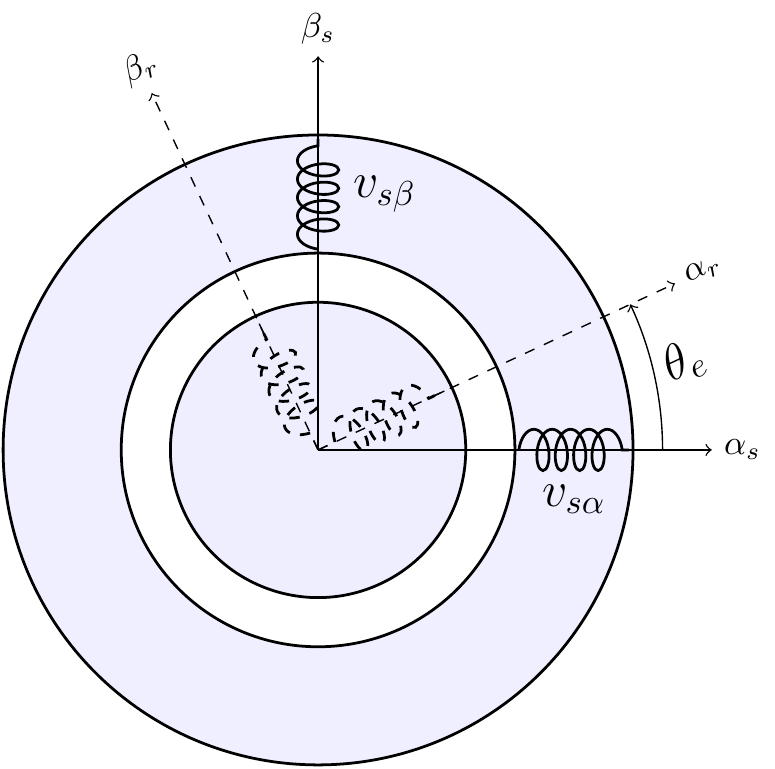} \caption{Schematic representation of the IM in the two-phase stationary reference frame}
	\label{im_schema}
\end{figure}

The machine model is expressed in the two-phase stationary reference frame $\alpha_s \beta_s$ (figure \ref{im_schema}), where the stator currents are directly measured. The state and input vectors are the following:
\begin{eqnarray*}
x = \begin{bmatrix}
{\mathcal{I}}_s^T & {\Psi}_r^T & \omega_e & T_r
\end{bmatrix}^T
~~~~;~~~~
u = \mathcal{V}_s
\end{eqnarray*}
${\mathcal{I}}_s$, ${\Psi}_r$ and $\mathcal{V}_s$ stand for the stator currents, rotor fluxes and stator voltages in $\alpha_s \beta_s$ coordinates:
\begin{equation*}
\mathcal{I}_s = \begin{bmatrix}
{i}_{s\alpha_s} & {i}_{s\beta_s}
\end{bmatrix}^T~~~~;~~~~
{\Psi}_r  = \begin{bmatrix}
{\psi}_{r\alpha_s} & {\psi}_{r\beta_s}
\end{bmatrix}^T~~~~;~~~~
{\mathcal{V}}_s = \begin{bmatrix}
{v}_{s\alpha_s} & {v}_{s\beta_s} 
\end{bmatrix}^T
\end{equation*}

Then, the state space model of the IM can be written as follows:
\begin{subequations}
\begin{eqnarray}
\frac{d \mathcal{I}_s}{dt} &=& - \frac{R_\sigma}{L_\sigma}  \mathcal{I}_s + \frac{k_r}{L_\sigma} \left(\frac{1}{\tau_r} \mathbf{I}_2 - \omega_e \mathbf{J}_2 \right) \Psi_r + \frac{1}{L_\sigma} \mathcal{V}_s \label{ss_IM:1}\\
\frac{d\Psi_{r}}{dt} &=& - \left(\frac{1}{\tau_r} \mathbf{I}_2 - \omega_e \mathbf{J}_2 \right) \Psi_r + \frac{M}{\tau_r} \mathcal{I}_s \label{ss_IM:2}\\
\frac{d \omega_e}{dt} &=& \frac{p^2}{J} k_r~\mathcal{I}_s^T \mathbf{J}_2 \Psi_r - \frac{p}{J}T_r \label{ss_IM:3}\\
\frac{d T_r}{dt} &=& 0 \label{ss_IM:4}
\end{eqnarray}
\label{ss_im}
\end{subequations}
where:
\begin{eqnarray*}
k_r &=& {M}/{L_r}~~~~;~~~~k_s = {M}/{L_s}~~~~;~~~~\tau_r = {L_r}/{R_r}\\
\sigma &=& 1 - k_r k_s ~~;~~ L_\sigma = \sigma L_s ~~;~~
R_\sigma = R_s + k_r^2 R_r
\end{eqnarray*}
$R$, $L$ and $M$ stand respectively for the resistance, inductance and mutual inductance. The indices $s$ and $r$ stand for stator and rotor parameters. $p$ is the number of pole pairs, $J$ is the inertia of the rotor with the associated load and $\omega_e$ is the electrical speed of the rotor. $\mathbf{I}_n$ is the $n \times n$ identity matrix, and $\mathbf{J}_2$ is the $\pi/2$ rotation matrix:
$$ \mathbf{J}_2 = \begin{bmatrix} 0 & -1 \\ 1 & 0 \end{bmatrix}$$

In order to facilitate the equations manipulation, the following change of variables is made:
\begin{eqnarray}
\widetilde{\mathcal{I}}_s &=&  L_\sigma \mathcal{I}_s\\
\widetilde{\Psi}_r &=& k_r \Psi_r
\end{eqnarray}

The system \eqref{ss_im} becomes:
\begin{subequations}
	\begin{eqnarray}
	\frac{d \mathcal{\widetilde{I}}_s}{dt} &=& \mathcal{V}_s + a \mathcal{\widetilde{I}}_s  + \gamma(t) {\widetilde{\Psi}_r} \label{di_s}\\
	\frac{d \widetilde{\Psi}_r}{dt} &=& - \gamma(t) {\widetilde{\Psi}_r} - \left(a-b\right) \mathcal{\widetilde{I}}_s \label{dpsi_r}\\
	\frac{d \omega_e}{dt} &=& \frac{c}{J}\mathcal{\widetilde{I}}_s^T\mathbf{J}_2 \widetilde{\Psi}_r - \frac{p}{J}T_r\\
	\frac{d T_r}{dt} &=& 0
	\end{eqnarray}
	\label{ss_im_2}
\end{subequations}
with
\begin{eqnarray*}
a = -\frac{R_\sigma}{L_\sigma} ~~;~~ b = -\frac{R_s}{L_\sigma} ~~;~~ c = \frac{p^2}{L_\sigma} ~~;~~ \gamma(t) = \left(\frac{1}{\tau_r} \mathbf{I}_2 - \omega_e \mathbf{J}_2\right) ~~;~~ \frac{d \gamma}{dt} = - \frac{d \omega_e}{dt} \mathbf{J}_2
\end{eqnarray*}

This new scaled model \eqref{ss_im_2} will be used for observability analysis below. It is a 6-${th}$ order model, the local observability study requires the evaluation of derivatives up to the 5-${th}$ order of the output. However, regarding the equations complexity, only the first and second order derivatives are evaluated for the IM. Higher order derivatives of the output are too lengthy, and very difficult to deal with. As the observability rank criterion provides sufficient conditions, information contained in the first and second order derivatives is rich enough to study the IM observability.



 
 \subsection{IM observability with speed measurement}
 Rotor fluxes and load torque of induction machines cannot be easily measured. They are often estimated using a state observer. In order to check the observability of these variables, the rotor speed is considered to be measured with the stator currents, and the stator voltages are considered to be known. The output vector is the following:
 \begin{equation*}
 y = \begin{bmatrix}
 \mathcal{\widetilde{I}}_s^T & \omega_e
 \end{bmatrix}^T
 \end{equation*}
 
The observability (sub-)matrix of the output and it first order derivative is:
 \begin{equation}
 \mathcal{O}_{y1}^{IM} = \frac{\partial}{\partial {x}}\begin{bmatrix}
 y \\ \dot{y}
 \end{bmatrix} = 
 \begin{bmatrix}
 1 & 0 & 0 & 0 & 0 & 0\\
 0 & 1 & 0 & 0 & 0 & 0\\
 0 & 0 & 0 & 0 & 1 & 0\\
 a & 0 & \frac{1}{\tau_r} & \omega_e & \widetilde{\psi}_{r \beta_s} & 0\\
 0 & a & -\omega_e & \frac{1}{\tau_r} & - \widetilde{\psi}_{r \alpha_s} & 0\\
 - \frac{c}{J} \widetilde{\psi}_{r \beta_s} & \frac{c}{J} \widetilde{\psi}_{r \alpha_s} & \frac{c}{J} \widetilde{i}_{s \beta_s} & -\frac{c}{J} \widetilde{i}_{s \alpha_s} & 0 & -\frac{p}{J}
 \end{bmatrix}
 \end{equation}
 Its determinant is the following:
 \begin{eqnarray}
 \Delta_{y1}^{IM} = -\frac{p}{J}\left(\omega_e^2 + \frac{1}{\tau_r^2}\right)
 \label{cond_obsv_im_w}
 \end{eqnarray}
We can then conclude that, for every stator voltages, currents and rotor speed and fluxes, the IM is observable if the rotor speed and stator currents are measured. There is no need to evaluate higher order derivatives.
 
\subsection{IM observability without speed measurement} 
In the context of sensorless control, only stator currents are measured, the output vector reads:
\begin{eqnarray*}
y = \widetilde{\mathcal{I}}_s
\end{eqnarray*}
The first order output derivative is:
\begin{eqnarray*}
\dot{y} = \frac{d\widetilde{\mathcal{I}}_s}{dt} = \mathcal{V}_s + a \mathcal{\widetilde{I}}_s  + \gamma(t) {\widetilde{\Psi}_r}
\end{eqnarray*}
Adding \eqref{di_s} and \eqref{dpsi_r} gives:
\begin{eqnarray*}
\frac{d\mathcal{\widetilde{I}}_s}{dt} + \frac{d\widetilde{\Psi}_r}{dt} = \mathcal{V}_s + b \mathcal{\widetilde{I}}_s
\end{eqnarray*}
then:
\begin{eqnarray*}
\frac{d\widetilde{\Psi}_r}{dt} = \mathcal{V}_s + b \mathcal{\widetilde{I}}_s  - \frac{d\mathcal{\widetilde{I}}_s}{dt}
\end{eqnarray*}

The second order derivative of the output can be then written as:
\begin{eqnarray*}
\frac{d^2\mathcal{\widetilde{I}}_s}{dt^2} &=& \frac{d\mathcal{\widetilde{V}}_s}{dt} + a \frac{d\mathcal{\widetilde{I}}_s}{dt} + \gamma(t) \frac{d\widetilde{\Psi}_r}{dt} + \frac{d\gamma}{dt}\widetilde{\Psi}_r\\
&=& \frac{d\mathcal{V}_s}{dt} + a \frac{d\mathcal{\widetilde{I}}_s}{dt} + \gamma(t) \mathcal{V}_s + \gamma(t) b \mathcal{\widetilde{I}}_s - \gamma(t) \frac{d\mathcal{\widetilde{I}}_s}{dt} + \frac{d\gamma}{dt} \widetilde{\Psi}_r \\
&=& \frac{d\mathcal{V}_s}{dt} + \gamma(t) \mathcal{V}_s + (a \mathbf{I}_2 - \gamma(t)) \frac{d \widetilde{\mathcal{I}}_s}{dt} + \gamma(t) b \mathcal{\widetilde{I}}_s + \frac{d\gamma}{dt} \widetilde{\Psi}_r
\end{eqnarray*}


The observability study is done using the scaled output and its derivatives:
\begin{eqnarray*}
h(x) &=& \mathcal{\widetilde{I}}_s\\
\mathcal{L}_fh(x) &=& \mathcal{V}_s + a \mathcal{\widetilde{I}}_s  + \gamma(t) {\widetilde{\Psi}_r}\\
\mathcal{L}^2_fh(x) &=& \frac{d\mathcal{V}_s}{dt} + a \mathcal{V}_s + \left(a^2 \mathbf{I}_2 - (a-b) \gamma(t) \right) \widetilde{\mathcal{I}}_s + \left(\frac{d\gamma}{dt} + a \gamma(t) - \gamma(t)^2\right) \widetilde{\Psi}_r 
\end{eqnarray*}

The IM observability (sub-)matrix, evaluated for these derivatives, can be written as:
\begin{eqnarray}
\mathcal{O}_{y2}^{IM} &=& \begin{bmatrix} 
1 & 0 & 0 & 0 & 0 & 0\\
0 & 1 & 0 & 0 & 0 & 0\\
a & 0 & \frac{1}{\tau_r} & \omega_e &  \widetilde{\psi}_{r \beta_s} & 0\\
0 & a & -\omega_e & \frac{1}{\tau_r} & -\widetilde{\psi}_{r \alpha_s} & 0\\  
d_{11} & d_{12} & e_{11} & e_{12} & f_{11} & f_{12}\\
d_{21} & d_{22} & e_{21} & e_{22} & f_{21} & f_{22} 
\end{bmatrix}
\label{IM_obsv_matrix}
\end{eqnarray}
with:
\begin{eqnarray*}
d_{11} &=& a^2 - \frac{a-b}{\tau_r} - \frac{c}{J} \widetilde{\psi}_{r \beta_s}^2 ~~~~~~;~~~~~~ d_{12} = -(a-b) \omega_e + \frac{c}{J} \widetilde{\psi}_{r \alpha_s} \widetilde{\psi}_{r \beta_s}\\
d_{21} &=& (a-b) \omega_e + \frac{c}{J} \widetilde{\psi}_{r \alpha_s} \widetilde{\psi}_{r \beta_s} ~~;~~~~~~ d_{22} = a^2 - \frac{a-b}{\tau_r} - \frac{c}{J} \widetilde{\psi}_{r \alpha_s}^2
\end{eqnarray*}
\begin{eqnarray*}
e_{11} &=& \frac{a}{\tau_r} - \frac{1}{\tau_r^2} + \omega_e^2 + \frac{c}{J} \widetilde{i}_{s \beta_s} \widetilde{\psi}_{r \beta_s} ~~~~~~~~~;~~~~~~ e_{12} = a \omega_e - 2 \frac{\omega_e}{\tau_r} + \frac{d\omega_e}{dt} - \frac{c}{J} \widetilde{i}_{s \alpha_s} \widetilde{\psi}_{r \beta_s}\\
e_{21} &=&  - a \omega_e + 2 \frac{\omega_e}{\tau_r} - \frac{d\omega_e}{dt} - \frac{c}{J} \widetilde{i}_{s \beta_s} \widetilde{\psi}_{r \alpha_s} ~;~~~~~~ e_{22} = \frac{a}{\tau_r} - \frac{1}{\tau_r^2} + \omega_e^2 + \frac{c}{J} \widetilde{i}_{s \alpha_s} \widetilde{\psi}_{r \alpha_s}
\end{eqnarray*}
\begin{eqnarray*}
f_{11} &=& 2 \omega_e \widetilde{\psi}_{r \alpha_s} - (a-b) \widetilde{i}_{s \beta_s} + \left(a - \frac{2}{\tau_r}\right) \widetilde{\psi}_{r \beta_s}  ~~~~~~;~~~~~~
f_{12} = -\frac{p}{J} \widetilde{\psi}_{r \beta_s}\\
f_{21} &=& 2 \omega_e \widetilde{\psi}_{r \beta_s} + (a-b) \widetilde{i}_{s \alpha_s} - \left(a + \frac{2}{\tau_r}\right) \widetilde{\psi}_{r \alpha_s} ~~~~~~;~~~~~~ f_{22} = \frac{p}{J} \widetilde{\psi}_{r \alpha_s}
\end{eqnarray*}

The matrix $\mathcal{O}_{y2}^{IM}$ \eqref{IM_obsv_matrix} is a $6\times6$ matrix, its determinant is the following:

\begin{eqnarray}
\Delta_{y2}^{IM} = \frac{p}{J} \frac{k_r^2}{\tau_r^2} \left[\tau_r \dot{\omega}_e\left({\psi}_{r \alpha_s}^2 + {\psi}_{r \beta_s}^2\right) - \left(1 + \tau_r^2 \omega_e^2\right) \left( \frac{d {\psi}_{r \alpha_s}}{dt} {\psi}_{r \beta_s} - \frac{d {\psi}_{r \beta_s}}{dt} {\psi}_{r \alpha_s}\right)\right]
\label{delta_im}
\end{eqnarray}
It is sufficient that the value of $\Delta_{y2}^{IM}$ is non-zero to guarantee the IM local observability. The expression \eqref{delta_im} suggests observability problems in the following two particular situations:
\begin{itemize}
\item $\dot{\omega}_e = 0$: For constant rotor speed, to ensure observability, it is sufficient that the rotor fluxes are neither constant nor linearly dependent, i.e. $\psi_{r\alpha_s} \neq k.\psi_{r\beta_s}$, for all constant $k$.
\item $\dot{\psi}_{r \alpha_s} = \dot{\psi}_{r \beta_s} =0$: If the rotor fluxes are slowly varying (constant), which usually corresponds to low-frequency input voltages, observability can be guaranteed by varying the rotor speed.
\end{itemize}  

Moreover, $\Delta_{y2}^{IM}$ can be expressed as follows:
\begin{eqnarray*}
\Delta_{y2}^{IM} = \frac{p}{J} \frac{k_r^2}{\tau_r^2} \left[ \tau_r \dot{\omega}_e \Psi_r^T \Psi_r - \left(1 + \tau_r^2 \omega_e^2 \right) \left({\Psi}_r^T \mathbf{J}_2 \dot{\Psi}_r \right)\right] 
\end{eqnarray*}
Then, the condition $\Delta_{y2}^{IM} = 0$ implies:
\begin{eqnarray*}
\frac{\tau_r \dot{\omega}_e}{1 + \tau_r^2 \omega_e^2} = \left(\Psi_r^T \Psi_r \right)^{-1} {\Psi}_r^T \mathbf{J}_2 \dot{\Psi}_r
\end{eqnarray*}
which is equivalent to the following:
\begin{eqnarray*}
\frac{d}{dt} \arctan(\tau_r \omega_e)  + \frac{d}{dt} \arctan\left(\frac{\psi_{r \beta}}{\psi_{r \alpha}}\right) = 0
\end{eqnarray*}
Furthermore, we know that:
\begin{eqnarray*}
\frac{d}{dt} \arctan\left(\frac{\psi_{r \beta}}{\psi_{r \alpha}}\right) = \frac{d}{dt} \left(\angle \Psi_r\right) = \omega_s
\end{eqnarray*}
where $\angle \Psi_r$ is the angular position of the rotor flux vector, and $\omega_s$ is the angular frequency of the rotating magnetic field. The observability condition of the IM can then be formulated as:
\begin{eqnarray}
\frac{d}{dt} \arctan(\tau_r \omega_e)  + \omega_s \neq 0
\label{cond_obsv_im}
\end{eqnarray}

\begin{figure}[!htb]
	\centering
\includegraphics[scale = 1.2]{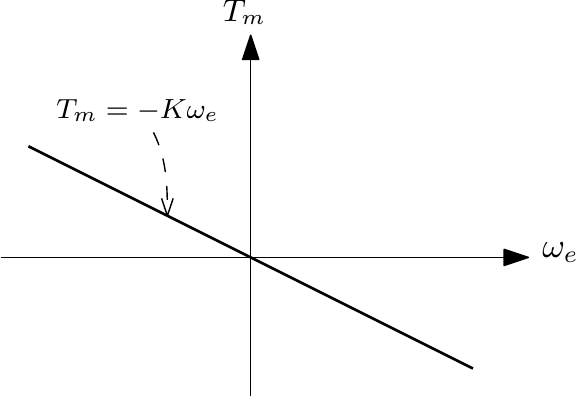}
	\caption{Unobservability line of the IM in the ($\omega_e,T_m$) plane. $K = p \frac{\psi_{r d_s}^2}{R_r}$}
	\label{obsv_line_im}
\end{figure}

The most important critical situation that can occur in practice is the low-frequency input voltage ($\omega_s = 0$). Recall the well-known relation in IM:
\begin{eqnarray*}
\omega_s = \omega_e + \omega_r
\end{eqnarray*}
where $\omega_r$ is the angular frequency of the rotor currents. It can be written as:
\begin{eqnarray*}
\omega_r = \frac{R_r}{p} \frac{T_m}{\psi_{r d_s}^2}
\end{eqnarray*}
where $T_m$ is the motor torque and $\psi_{rd_s}$ is the direct rotor flux on the rotor field-oriented reference frame.
Therefore, the condition $\omega_s = 0$ is equivalent to:
\begin{eqnarray}
\omega_e = - \omega_r ~~\Leftrightarrow ~~ \omega_e = - \frac{R_r}{p} \frac{T_m}{\psi_{r d_s}^2}
\end{eqnarray}
This equation describes a line in the ($\omega_r$,$T_m$) plane (figure \ref{obsv_line_im}); it is located in the second and fourth quadrants, that correspond to the machine operating in generator mode. By abuse of notation, this line will be called "unobservability line", although one can only say that observability is not \emph{guaranteed} along this line.


\section{Simulation results}
In this section, numerical simulations using \emph{Matlab/Simulink} are carried out in order to examine the above observability conditions, for both the WRSM and the IM. An extended Kalman filter is designed for this purpose. As only the observability is concerned, the observer is operated in open loop to avoid system's stability issues. The simulation scenarios are chosen in a way to show the observer's behavior in the critical observability situations. The model parameters of the machines are given in the appendix \ref{appendix_a}.

\subsection{WRSM}
\begin{figure}[!hb]
	\centering
	\includegraphics[width = 0.85\linewidth]{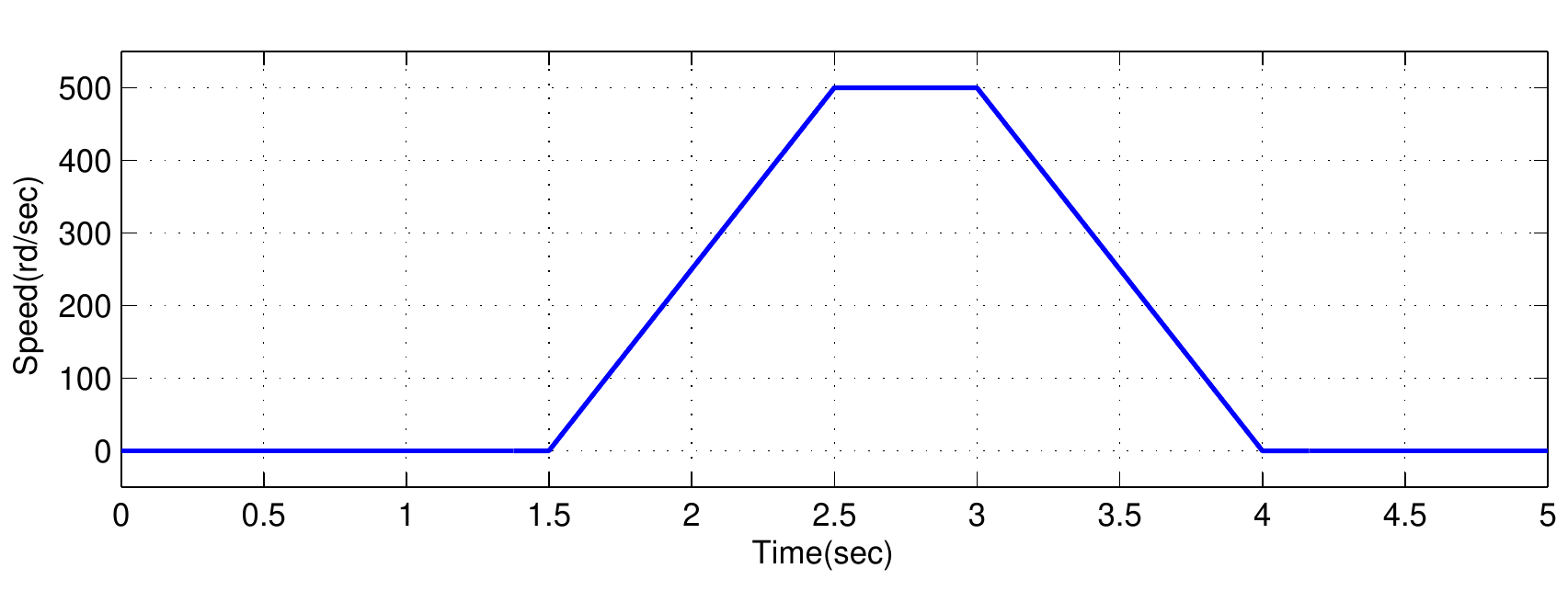}
	\caption{Rotor speed profile of the WRSM}
	\label{speed_WRSM}
\end{figure}

The following simulation scenario is considered for the WRSM: The machine's shaft is driven by an external mechanical system that forces the speed profile shown in the figure \ref{speed_WRSM}. The currents $i_d$ and $i_q$ are calculated using the real rotor position, and fed back to the current controller. Proportional-integral (PI) controllers are designed for the currents to fit with the following set-points:
\begin{eqnarray*}
i_d^{\#} = 2~A ~~~~~~;~~~~~~ i_q^{\#} = 15~A ~~~~~~;~~~~~~ i_f^{\#} = 4~A
\end{eqnarray*}  

The EKF tuning matrices are taken as follows:
\begin{eqnarray*}
\mathbf{Q} = \begin{bmatrix}
1 & 0 & 0 & 0 & 0\\
0 & 1 & 0 & 0 & 0\\
0 & 0 & 1 & 0 & 0\\
0 & 0 & 0 & 200 & 0\\
0 & 0 & 0 & 0 & 5
\end{bmatrix} ~~~~~~;~~~~~~
\mathbf{R} = \begin{bmatrix}
1 & 0 & 0 \\
0 & 1 & 0 \\
0 & 0 & 1
\end{bmatrix}
\end{eqnarray*}

The aim is to examine the observer behavior at standstill, for both fixed and moving observability vector situations. For this reason, a high frequency current is superimposed to the rotor current $i_f$ during the time intervals [$1$ $s$; $1.5$ $s$] and [$4.5$ $s$; $5$ $s$]; the HF signal is added to the $i_f$ set-point as follows:
\begin{eqnarray*}
i_f^{\#} = i_{f0} + i_{f HF} = I_{f0} + I_{f HF} \sin \omega_{HF}t
\end{eqnarray*}
where 
\begin{eqnarray*}
I_{f0} = 4~A~~~~;~~~~I_{f HF} = 0.5~A ~~~~;~~~~ \omega_{HF} = 2 \pi. 10^3 ~rd/s
\end{eqnarray*}


\begin{figure}[!hb]
	\centering
	\includegraphics[width = 0.85\linewidth]{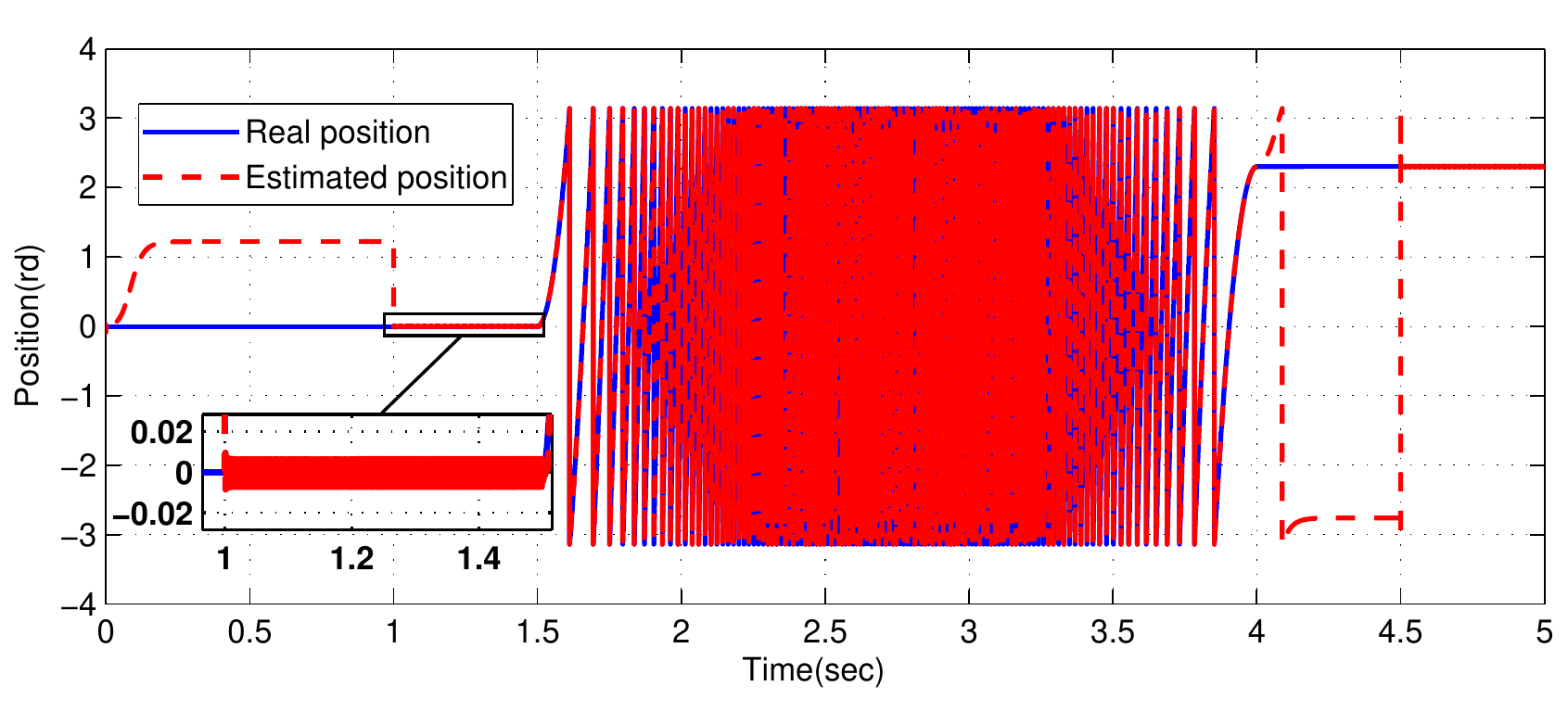}
	\caption{Rotor position estimation of the WRSM}
	\label{position_estimation_hf_WRSM}
\end{figure}

\begin{figure}[!hb]
	\centering
	\includegraphics[width = 0.85\linewidth]{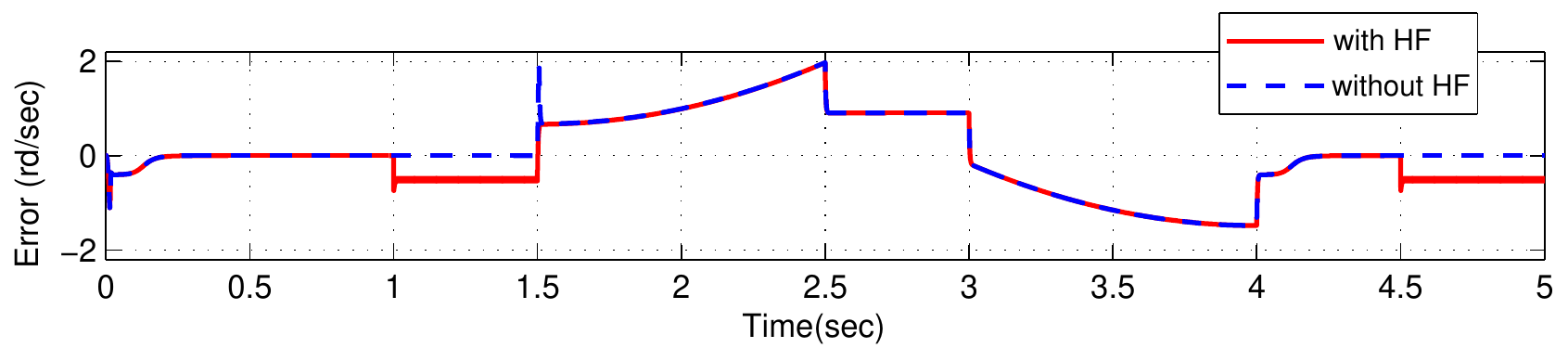}
	\caption{Rotor speed estimation error of the WRSM}
	\label{speed_estimation_WRSM}
\end{figure}

\begin{figure}[!hb]
	\centering
	\includegraphics[width = 0.85\linewidth]{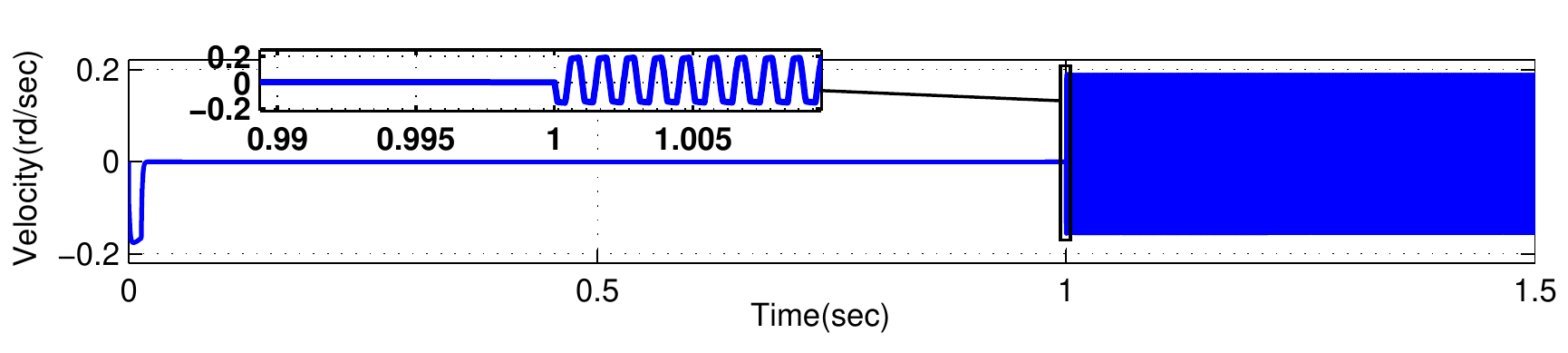}
	\caption{Angular velocity of the observability vector of the WRSM}
	\label{obsv_vector_WRSM}
\end{figure}

The position estimation is shown on the figure \ref{position_estimation_hf_WRSM}: at standstill, the EKF does not converge to the correct position until the HF current is injected. For nonzero speeds, there is no position estimation problem. The speed estimation error is shown on the figure \ref{speed_estimation_WRSM}: the error slightly increases with the HF injection, but it remains reasonable. The angular velocity of the observability vector is shown on the figure \ref{obsv_vector_WRSM} for the time interval [$0$ $s$; $1.5$ $s$]. It can be seen that this velocity is different from zero when the HF signal is injected, this proves that the simulation results are consistent with the observability analysis results of Section 3; the observer converges to the real position when the observability vector angular velocity is different from the rotor speed (in this case zero speed).

\subsection{IM}
\begin{figure}[!hb]
\centering
	\includegraphics[width = 0.85\linewidth]{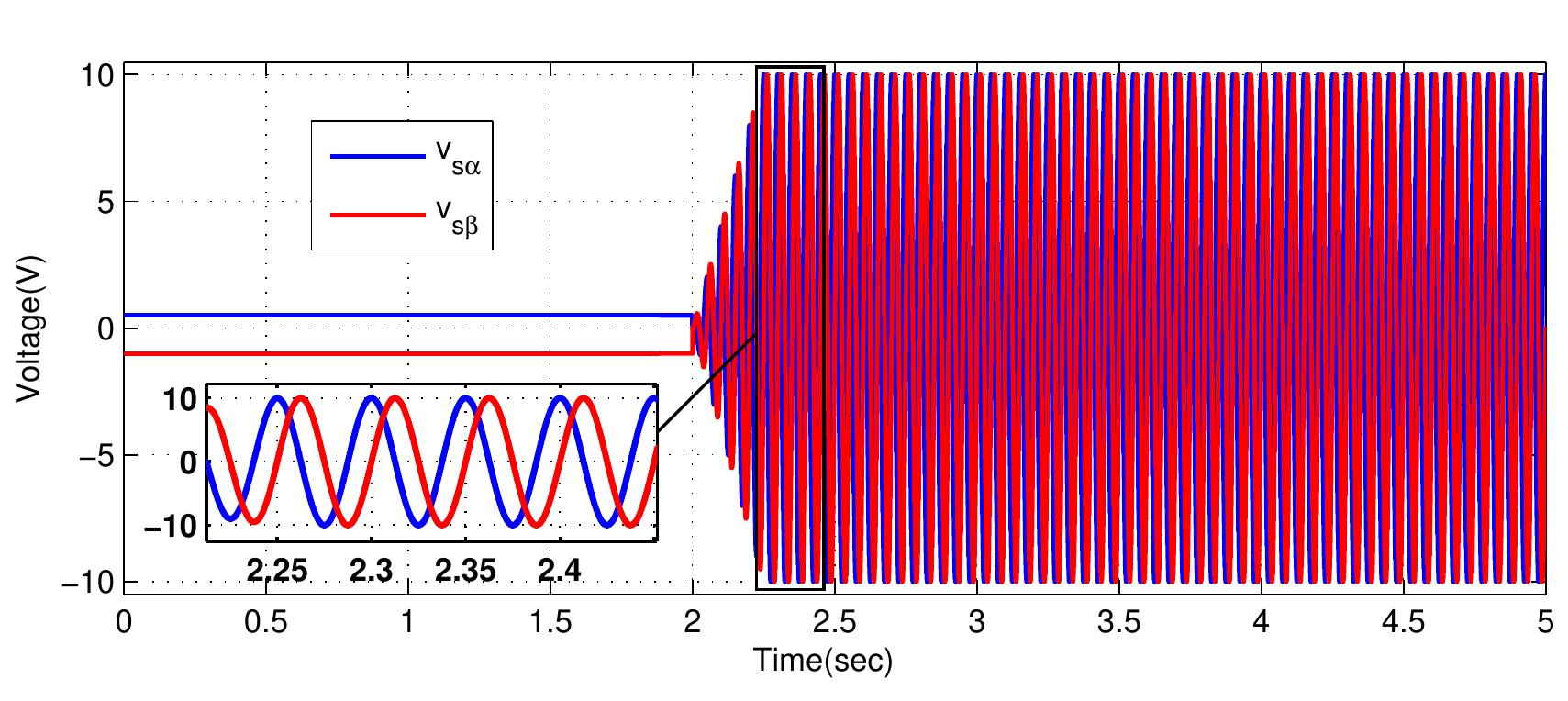}
	\caption{Input voltage of the IM}
	\label{v_sab_IM}
\end{figure}

\begin{figure}[!hb]
	\centering
	\includegraphics[width = 0.85\linewidth]{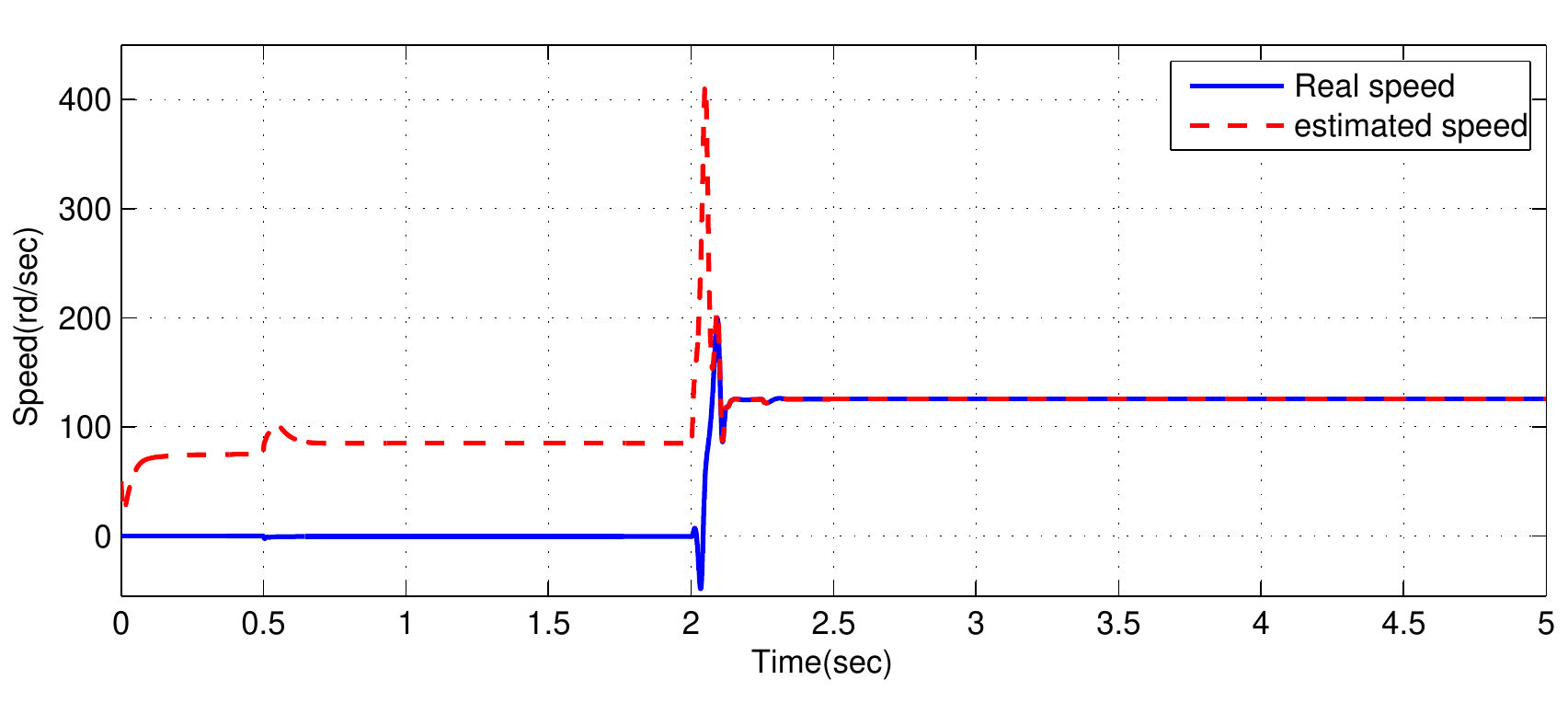}
	\caption{Rotor speed estimation of the IM}
	\label{speed_estimation_IM}
\end{figure}

\begin{figure}[!hb]
	\centering
	\includegraphics[width = 0.85\linewidth]{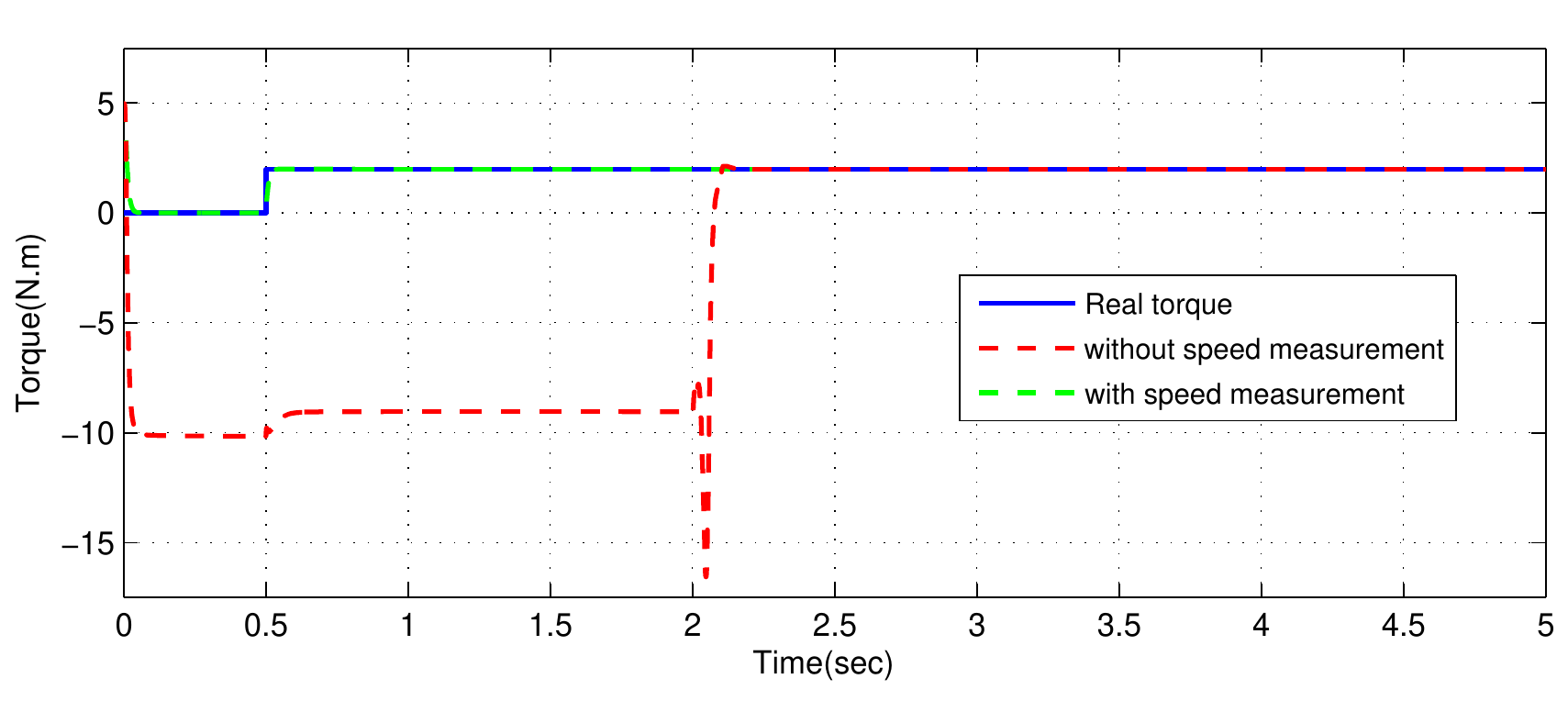}
	\caption{Load torque estimation of the IM}
	\label{torque_estimation_IM}
\end{figure}

\begin{figure}[!hb]
	\centering
	\includegraphics[width = 0.85\linewidth]{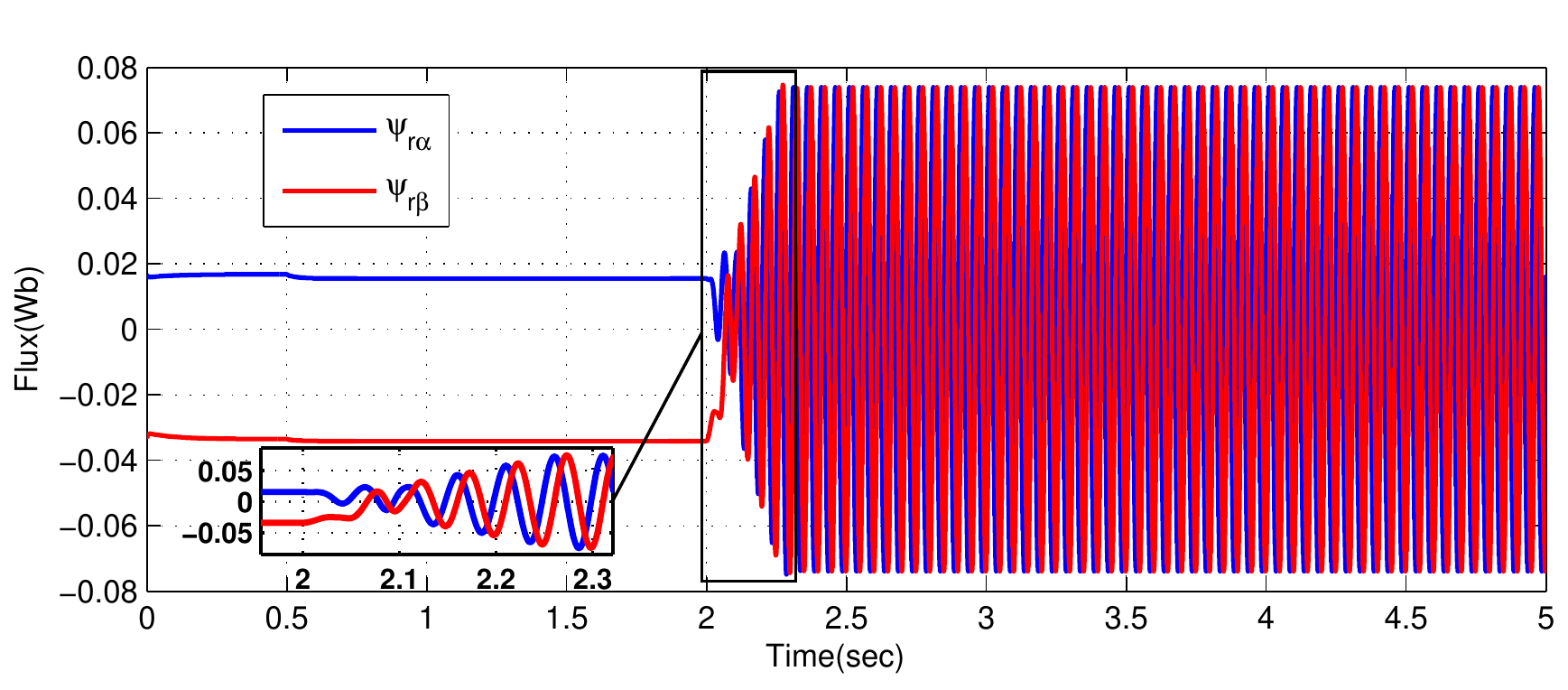}
	\caption{Rotor fluxes of the IM}
	\label{psi_rab_IM}
\end{figure}

The following simulation scenario is considered for the IM: the machine is operating in open loop. The applied stator voltages are shown in the figure \ref{v_sab_IM} in the two-phase stationary reference frame. A load torque profile, shown in the figure \ref{torque_estimation_IM}, in applied to the machine's rotor. The resulting rotor speed profile and rotor fluxes of the machine are presented, respectively, in figures \ref{speed_estimation_IM} and \ref{psi_rab_IM}. This scenario allows us to analyze the IM observability for zero angular frequency $\omega_s$. An EKF in designed to examine the IM observability in two cases: with and without rotor speed measurement. It is expected that the observer converges from the beginning when the speed is measured. The aim is to monitor the observer's behavior for zero-stator-frequency when the rotor speed is not measured.

In order to make the comparison easier, the same matrix $\mathbf{Q}$ is taken for both cases:
\begin{eqnarray*}
\mathbf{Q} = \begin{bmatrix}
100 & 0 & 0 & 0 & 0 & 0\\
0 & 100 & 0 & 0 & 0 & 0\\
0 & 0 & 0.1 & 0 & 0 & 0\\
0 & 0 & 0 & 0.1 & 0 & 0\\
0 & 0 & 0 & 0 & 10^{5} & 0\\
0 & 0 & 0 & 0 & 0 & 10^{4}
\end{bmatrix}
\end{eqnarray*}
and the matrix $\mathbf{R}$ used in the case where the rotor speed is measured is:
\begin{eqnarray*}
\mathbf{R} = 100~\mathbf{I}_3
\end{eqnarray*}
When the speed is not measured, the matrix $\mathbf{R}$ is the following:
\begin{eqnarray*}
\mathbf{R} = 100~\mathbf{I}_2
\end{eqnarray*}
The following initial conditions are introduced in the observer:
\begin{eqnarray*}
\hat{x}(0) &=& \begin{bmatrix}
0~A & 0~A & -0.02~Wb & -0.02~Wb & 50~rd/s & 5~N.m
\end{bmatrix}^T
\\
\mathbf{P}(0) &=& 0.1~\mathbf{I}_6
\end{eqnarray*}

\begin{figure}[!hb]
	\centering
\subfloat[$\psi_{r \alpha}$]{	\includegraphics[width = 0.85\linewidth]{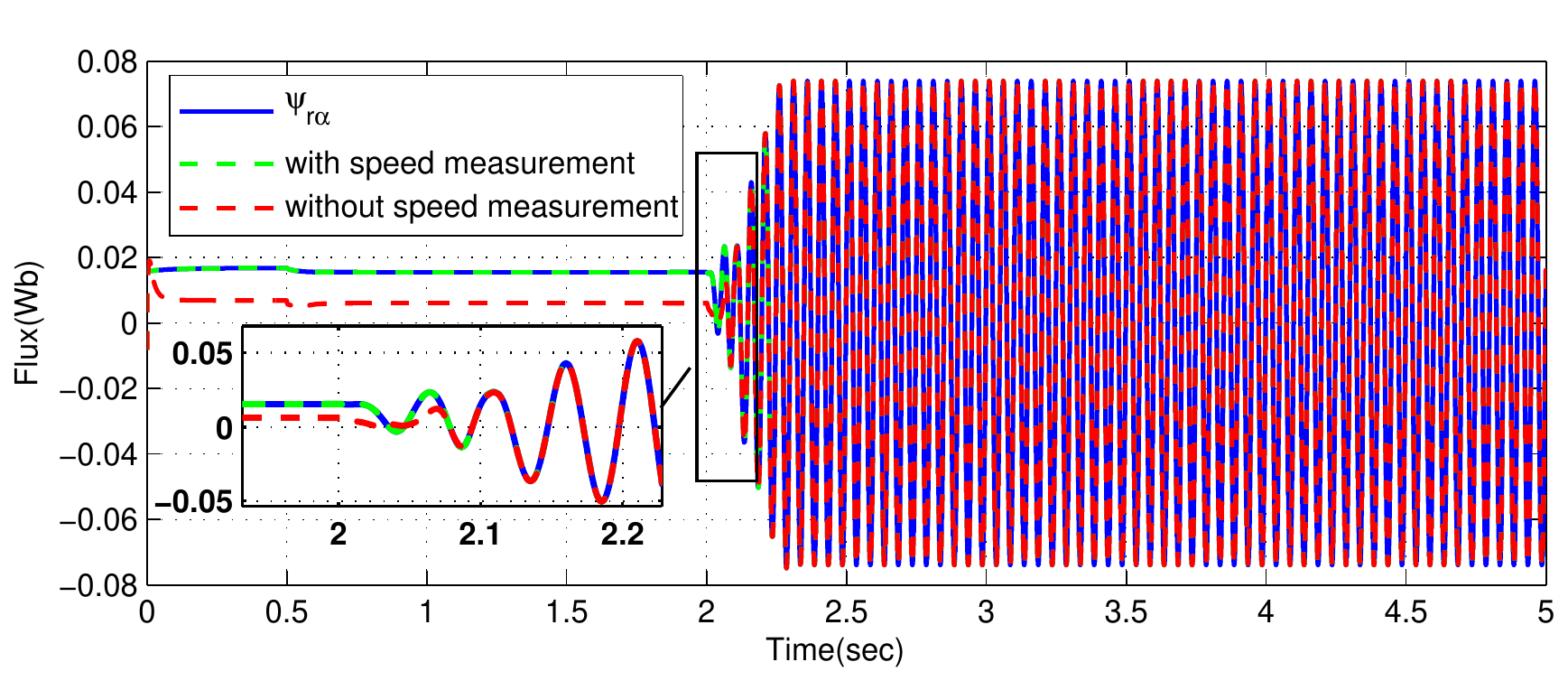}
	\label{psi_ra_est_IM}}

\subfloat[$\psi_{r \beta}$]{	\includegraphics[width = 0.85\linewidth]{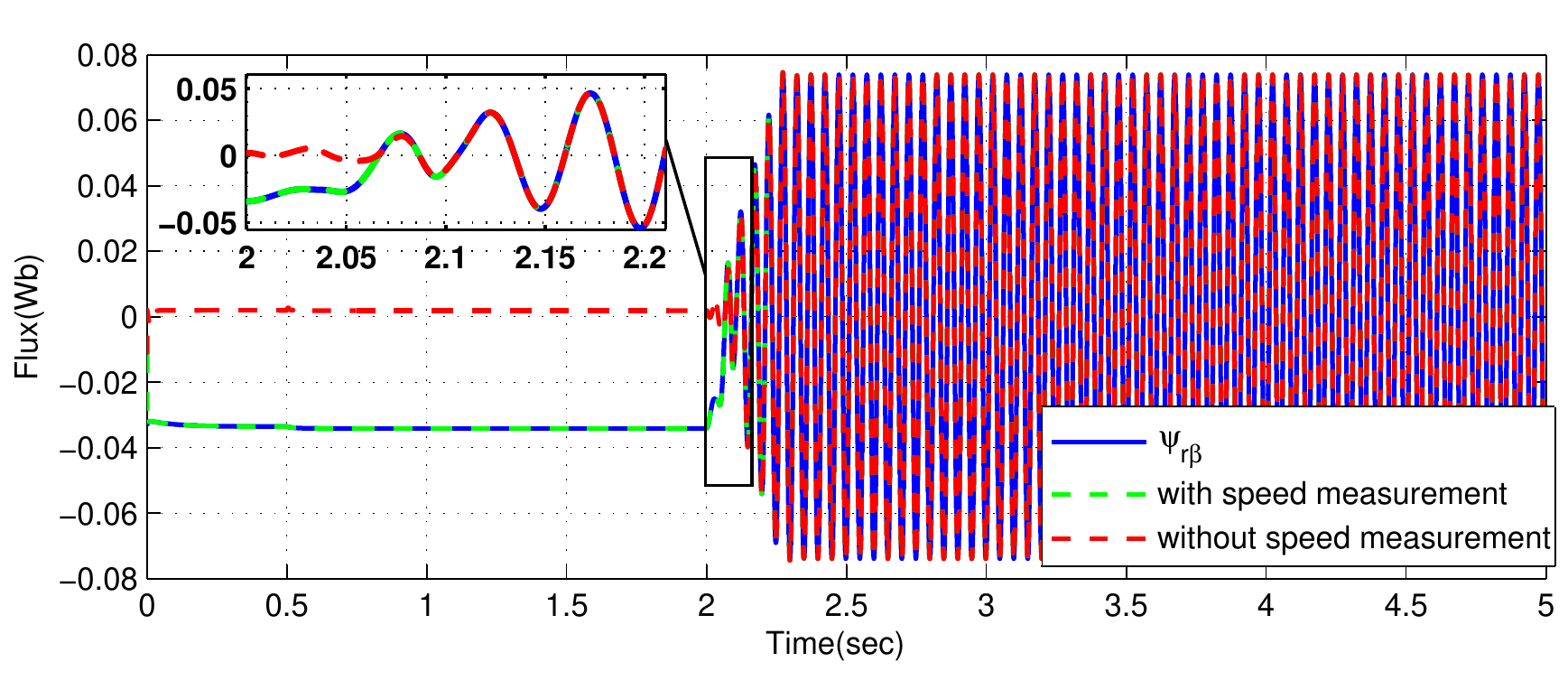}
	\label{psi_rb_est_IM}}
	\caption{Rotor flux estimation of the IM}
\end{figure}

The rotor fluxes estimation is shown in figure \ref{psi_ra_est_IM} and \ref{psi_rb_est_IM} for $\psi_{r\alpha}$ and $\psi_{r \beta}$ respectively. And the rotor speed and load torque estimation are shown in figure \ref{speed_estimation_IM} and \ref{torque_estimation_IM} respectively\footnote{Note that it is not important to show the speed estimation when it is measured.}. As expected, the observer converges for the whole scenario with rotor speed measurement. On the other hand, when the rotor speed is not measured, the fluxes, speed and torque are not correctly estimated if $\omega_s$ is null.

\section{Experimental results}
In this section, the observability conditions are checked on an experimental scenario. We recall that the observer is operated in open-loop, in order to distinguish observability issues from stability issues. Furthermore, the observer \emph{fine tuning} is beyond the scope of this paper. The machines parameters are the same as those used in simulation with some uncertainty; they are given in the appendix \ref{appendix_a}. The experimental data are taken from an experimental set-up of \emph{Renault} and are treated using the \emph{Matlab/Simulink} environment. 
 
 \subsection{WRSM}
\begin{figure}[!hb]
\centering
\includegraphics[width = 0.85\linewidth]{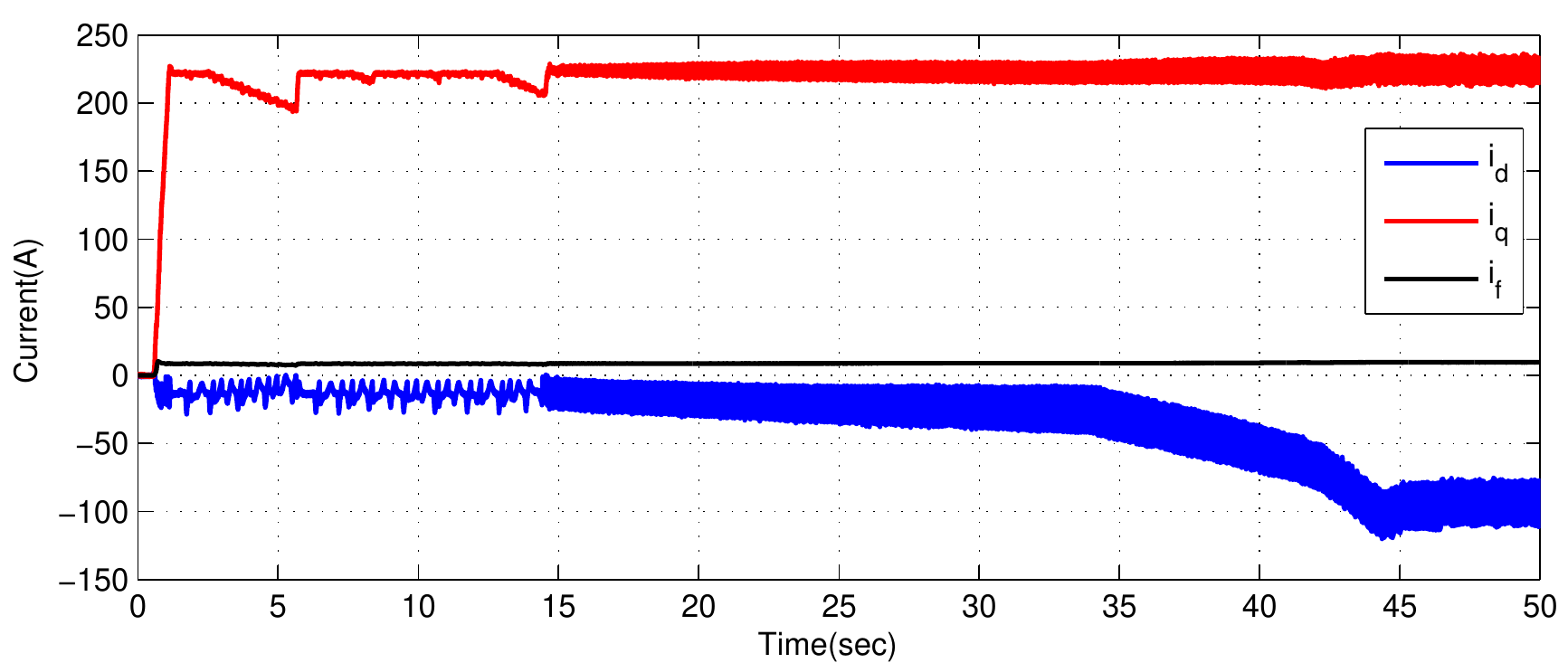}
	\caption{WRSM currents in the $dq$ coordinates}
	\label{i_dqf_exp_wrsm}
\end{figure}

\begin{figure}[!hb]
	\centering
	\includegraphics[width = 0.85\linewidth]{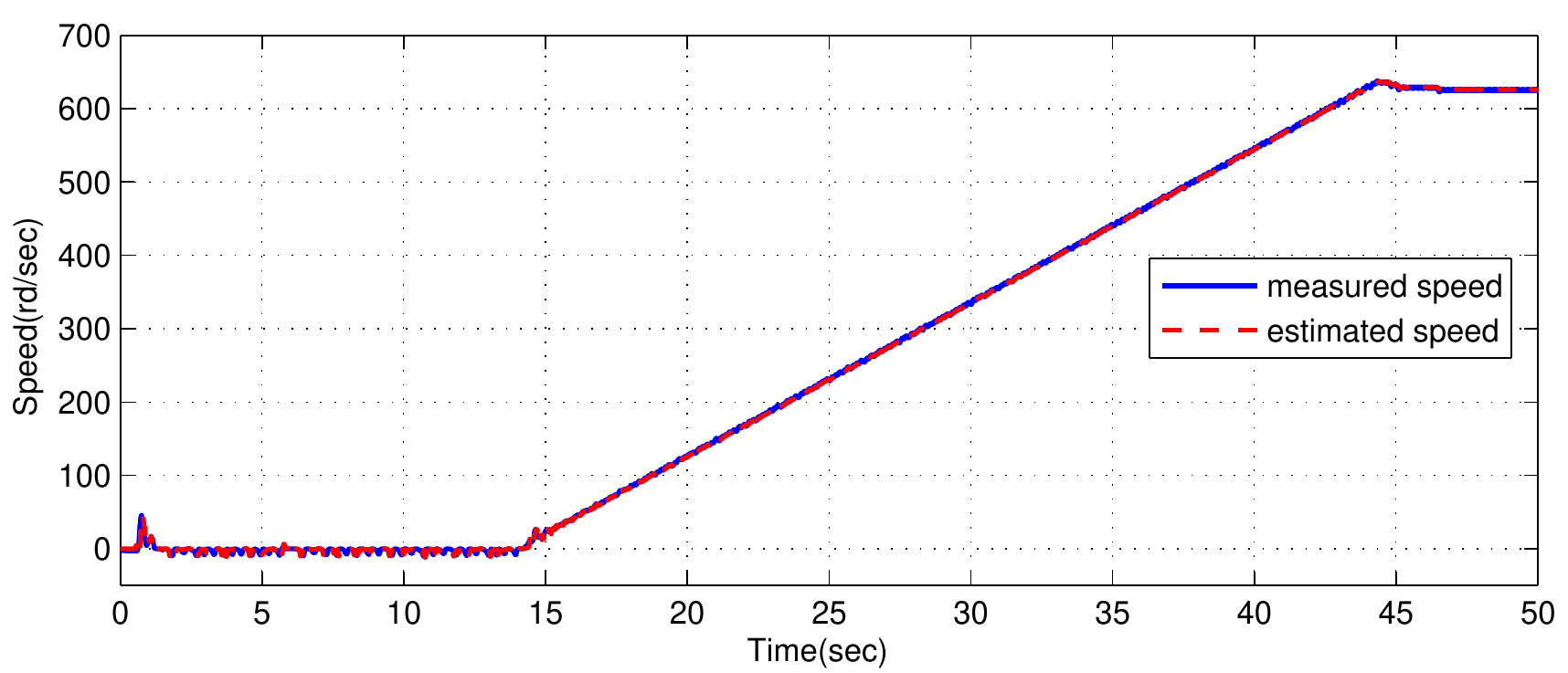}
	\caption{Rotor speed estimation of the WRSM}
	\label{omega_est_exp_wrsm}
\end{figure}

The experimental scenario for the WRSM is the following: the measured currents of the WRSM are given on the figure \ref{i_dqf_exp_wrsm} in the $dq$ reference frame, using the measured position. The speed profile used to drive the machine's shaft is illustrated in the figure \ref{omega_est_exp_wrsm}. This scenario enables us to check the observability conditions at standstill.

The Kalman tuning matrices are chosen as follows: 
\begin{eqnarray*}
\mathbf{Q} = \begin{bmatrix}
10^{10} & 0 & 0 & 0 & 0\\
0 & 10^{10} & 0 & 0 & 0\\
0 & 0 & 0.02 & 0 & 0\\
0 & 0 & 0 & 2 & 0\\
0 & 0 & 0 & 0 & 10^{-8}
\end{bmatrix} ~~~;~~~
\mathbf{R} = \begin{bmatrix}
10^{8} & 0 & 0 \\
0 & 10^{8} & 0 \\
0 & 0 & 1
\end{bmatrix}
\end{eqnarray*}
The sampling period is $T_s = 10^{-4}s$.

\begin{figure}[!hb]
	\centering
\includegraphics[width = 0.85\linewidth]{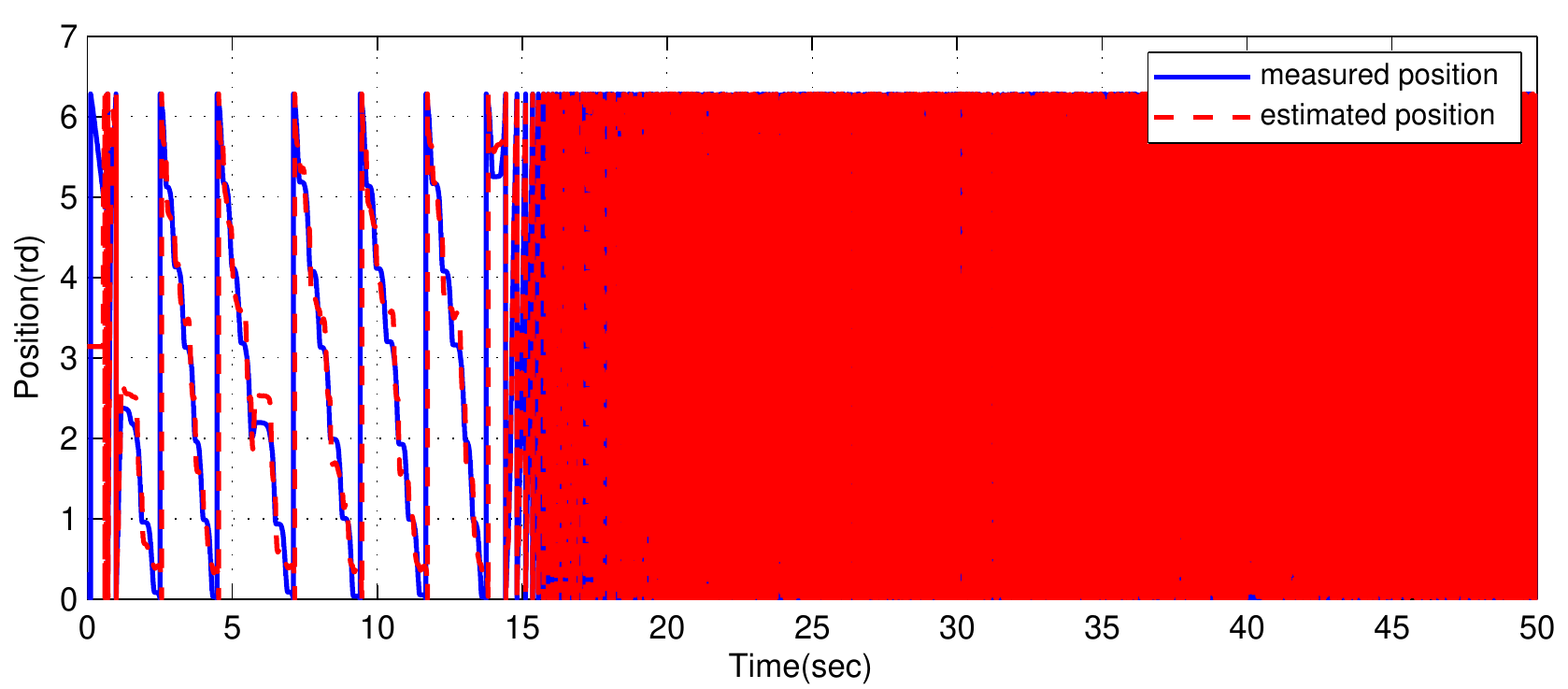}
	\caption{Rotor position estimation of the WRSM}
	\label{theta_exp_wrsm}
\end{figure}

\begin{figure}[!hb]
	\centering
\includegraphics[width = 0.85\linewidth]{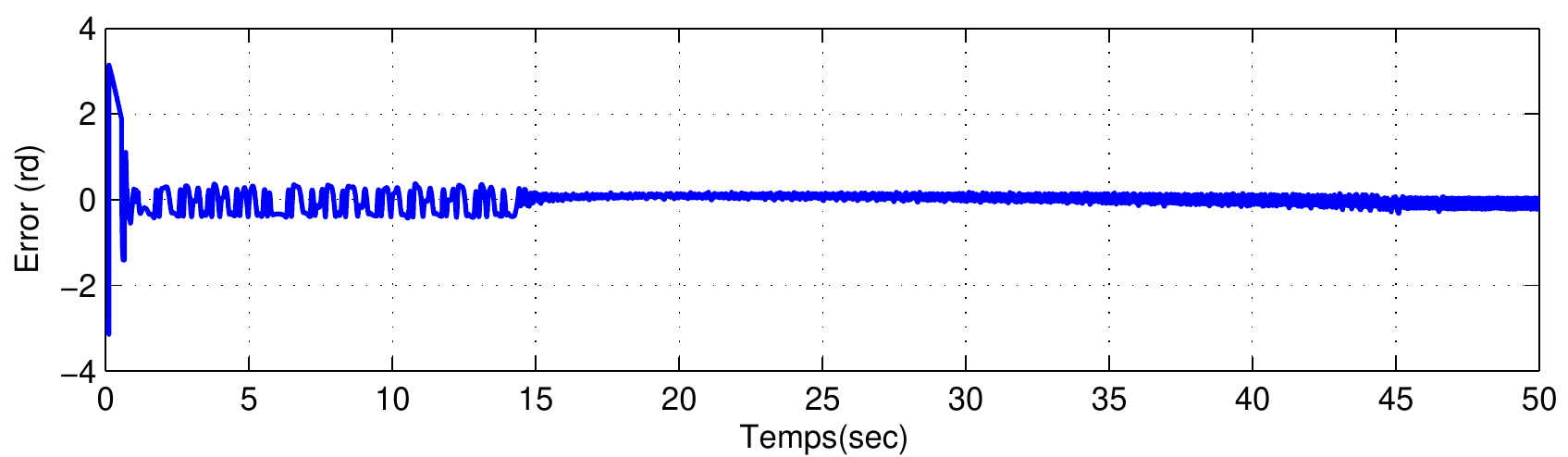}
	\caption{Rotor position estimation error of the WRSM}
	\label{theta_er_exp_wrsm}
\end{figure}

\begin{figure}[!hb]
	\centering
	\includegraphics[width = 0.85\linewidth]{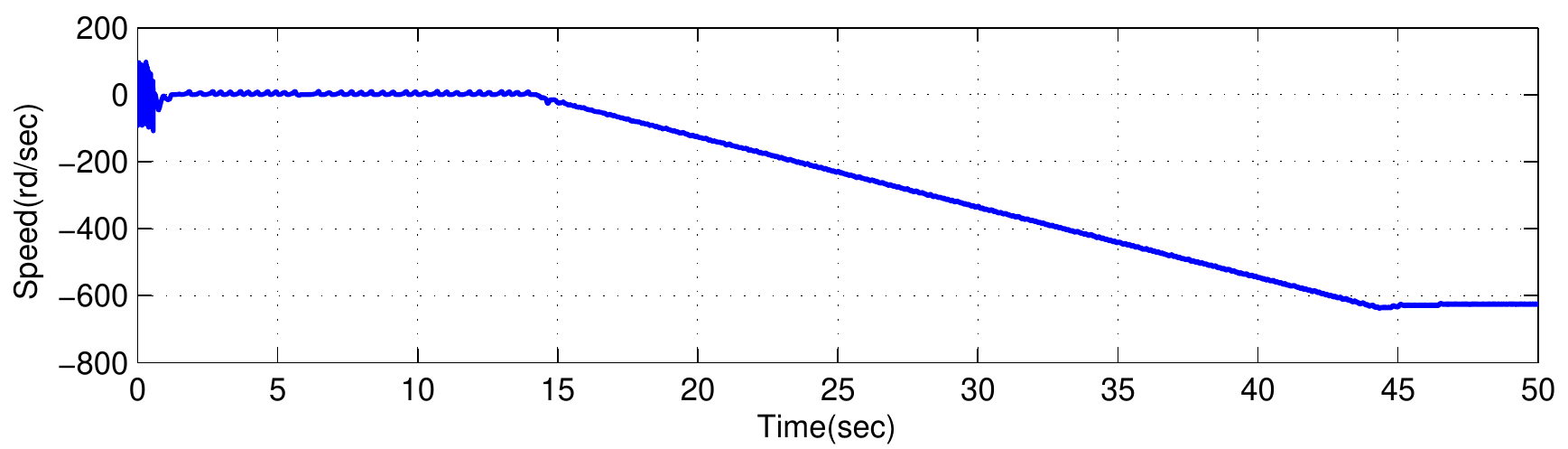}
	\caption{Difference between the WRSM rotor angular velocity and the observability vector one}
	\label{omega_o_exp_wrsm}
\end{figure}

\begin{figure}[!hb]
	\centering
\subfloat[Rotor speed estimation error of the WRSM]{\includegraphics[width = 0.85\linewidth]{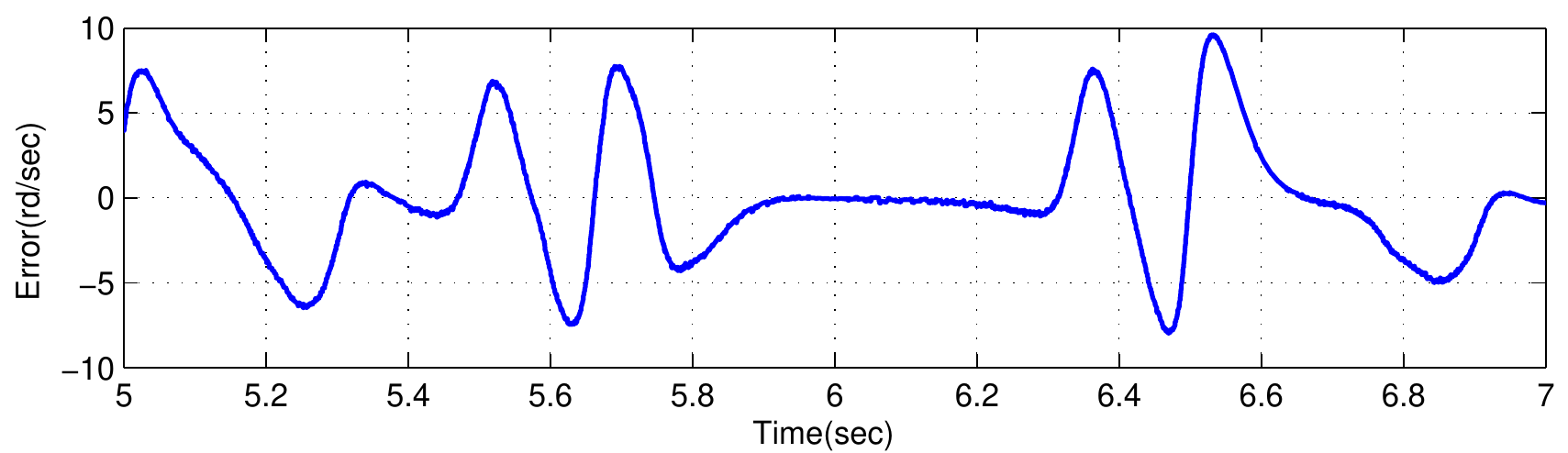}
	\label{omega_er_exp_6_wrsm}}

\subfloat[Rotor position estimation error of the WRSM]{	\includegraphics[width = 0.85\linewidth]{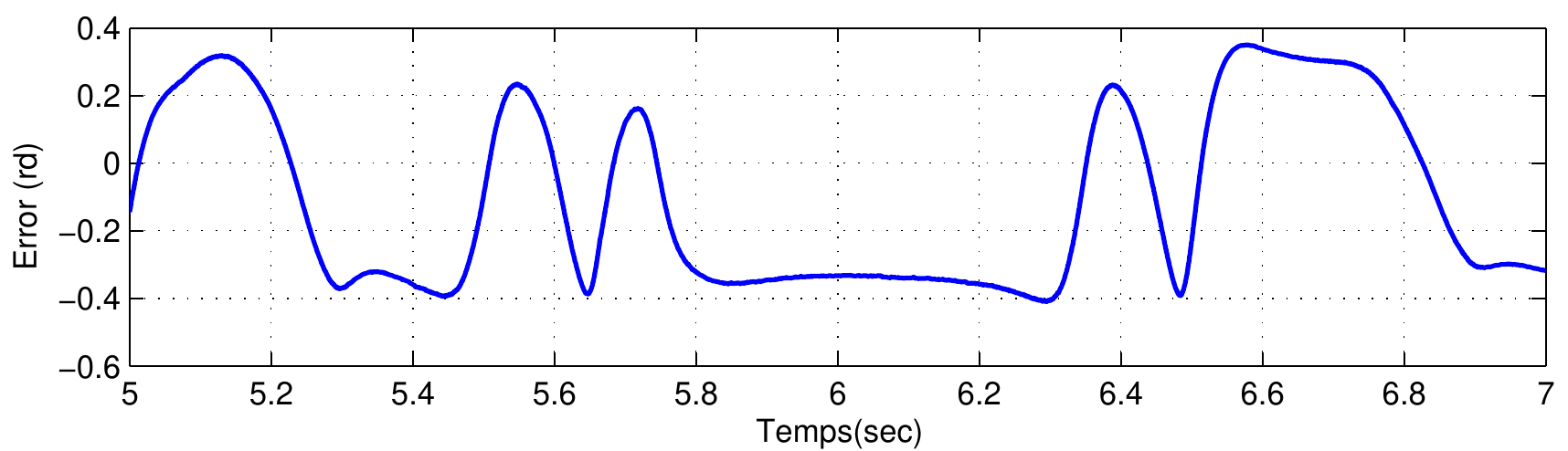}
	\label{theta_er_exp_6_wrsm}}

\subfloat[Difference between the WRSM rotor angular velocity and the observability vector one]{\includegraphics[width = 0.85\linewidth]{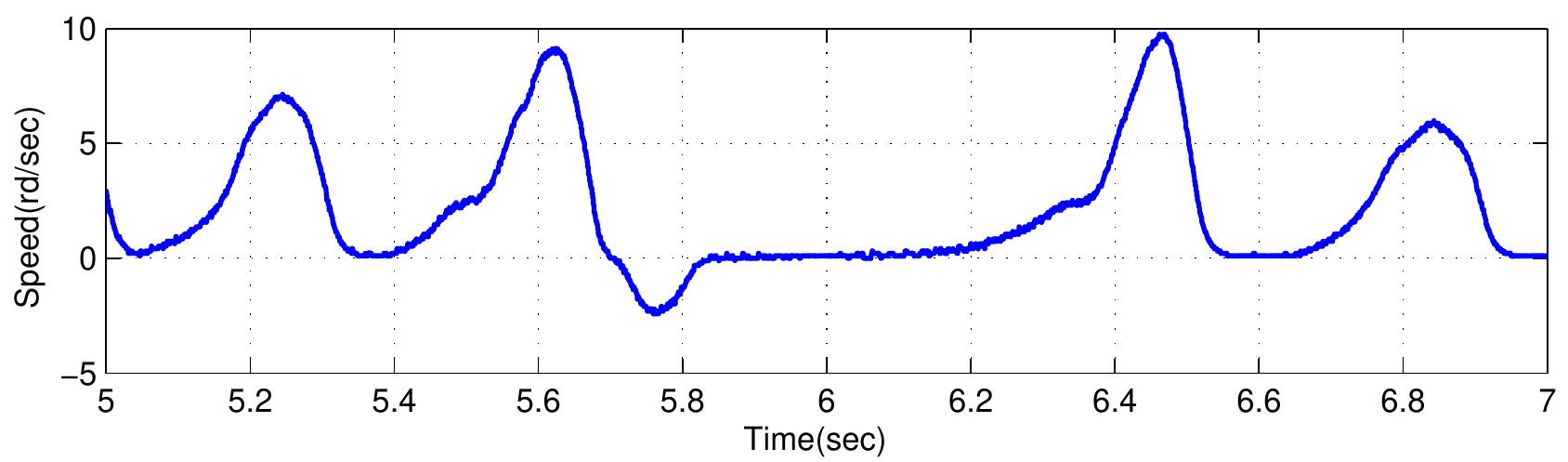}
	\label{omega_o_6_wrsm}}
    \caption{Estimation error on the time interval [$5$ $s$, $7$ $s$]}
\end{figure}

Figure \ref{omega_est_exp_wrsm} shows the measured speed against the estimated one. 
The position estimation is shown on the figure \ref{theta_exp_wrsm}, and the position estimation error on the figure \ref{theta_er_exp_wrsm}. The difference between angular velocities of the rotor and the observability vector, which is the image of the determinant of the observability matrix, is shown on the figure \ref{omega_o_exp_wrsm}. In order to analyze the standstill operating condition, the figures \ref{omega_er_exp_6_wrsm}, \ref{theta_er_exp_6_wrsm} and \ref{omega_o_6_wrsm} show, respectively, the rotor speed estimation error, the rotor position estimation error and the difference between the angular velocities of the rotor and the observability vector, for the time interval [$5~s$,$7~s$]: it can be seen that the position estimation error dynamics follows the relative angular velocity of the rotor and the observability vector, and that when the determinant is zero the position estimation error is different from zero (see instant $6~s$), whereas the speed estimation error is zero. This is consistent with our previous conclusions.

As a conclusion, the observability vector concept seems to be a useful tool for analyzing the performance of sensorless synchronous drives.

\subsection{IM}
\begin{figure}[!hb]
	\centering
	\includegraphics[width = 0.85\linewidth]{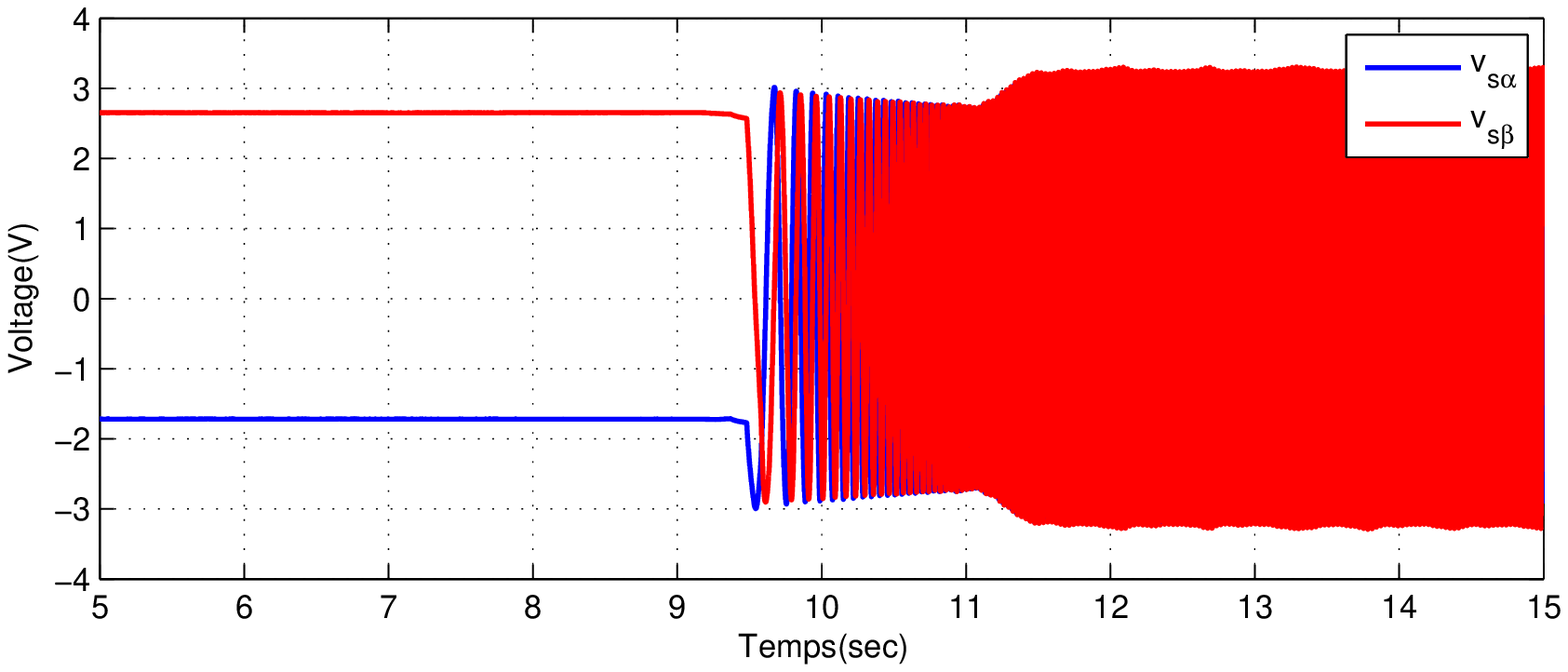}
	\caption{Stator voltages of the IM}
	\label{v_sab_exp_im_30}
\end{figure}

\begin{figure}[!hb]
	\centering
	\includegraphics[width = 0.85\linewidth]{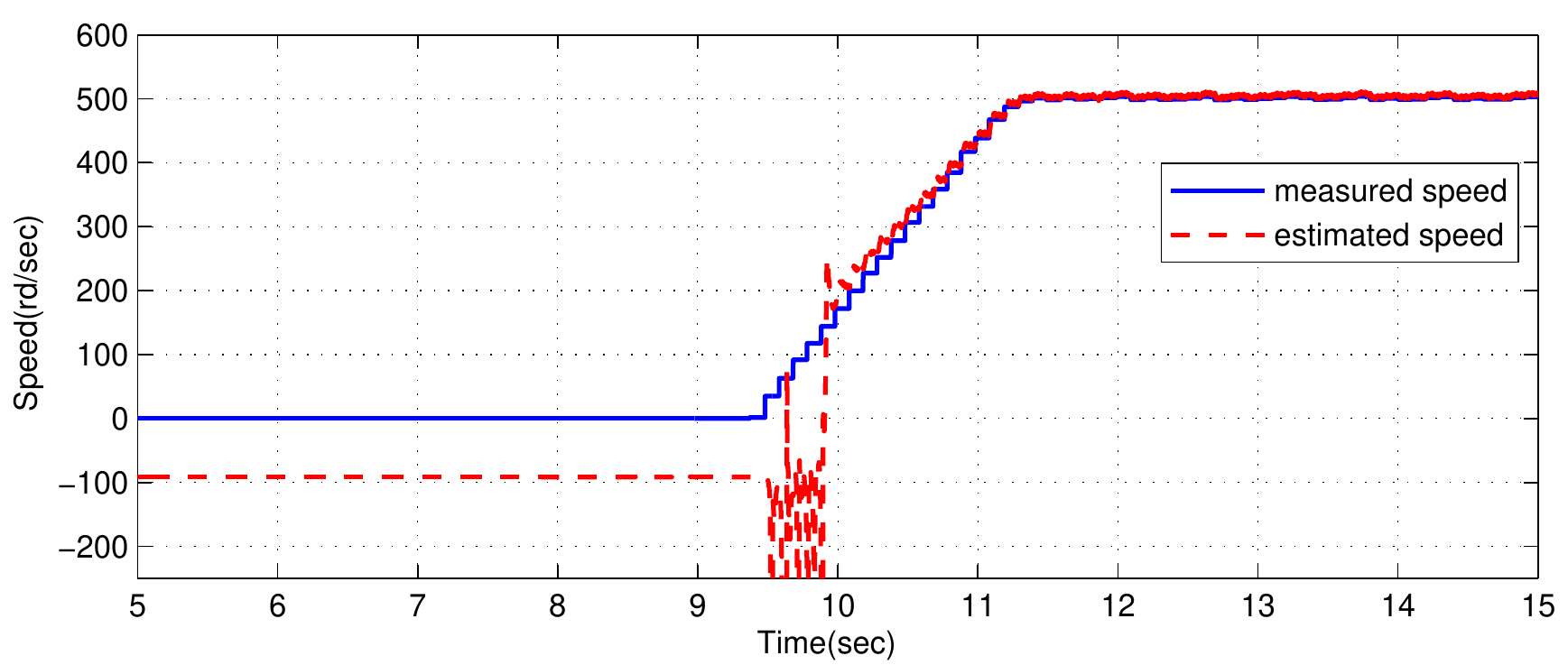}
	\caption{Rotor speed estimation}
	\label{omega_est_exp_im}
\end{figure}

The experimental scenario for the IM is the following: the stator voltages are shown on the  figure \ref{v_sab_exp_im_30} in the two-phase stationary reference frame, and the rotor speed profile is shown on the figure \ref{omega_est_exp_im}. In the sequel, the study is focused on the time interval [$5~s; 15~s$]. Two EKF filters are designed for observability study with and without speed measurement. The $\mathbf{Q}$ matrix is the same in both cases:

\begin{eqnarray*}
\mathbf{Q} = \begin{bmatrix}
10^5 & 0 & 0 & 0 & 0 & 0\\
0 & 10^5 & 0 & 0 & 0 & 0\\
0 & 0 & 10^{-3} & 0 & 0 & 0\\
0 & 0 & 0 & 10^{-3} & 0 & 0\\
0 & 0 & 0 & 0 & 10^{4} & 0\\
0 & 0 & 0 & 0 & 0 & 10^{4}
\end{bmatrix}
\end{eqnarray*}
When the speed is measured (output vector dimension equals $3$), the $\mathbf{R}$ matrix is tuned to:
\begin{eqnarray*}
\mathbf{R} = \begin{bmatrix}
10^{10} & 0 & 0\\
0 & 10^{10} & 0\\
0 & 0 & 10^5
\end{bmatrix}
\end{eqnarray*}
whereas without speed measurement, $\mathbf{R}$ is:
\begin{eqnarray*}
\mathbf{R} = 10^{10}~\mathbf{I}_2
\end{eqnarray*}
The sampling time is $T_s = 5 \times 10^{-5} s$.

\begin{figure}[!hb]
	\centering
\subfloat[$\psi_{r \alpha_s}$]{	\includegraphics[width = 0.85\linewidth]{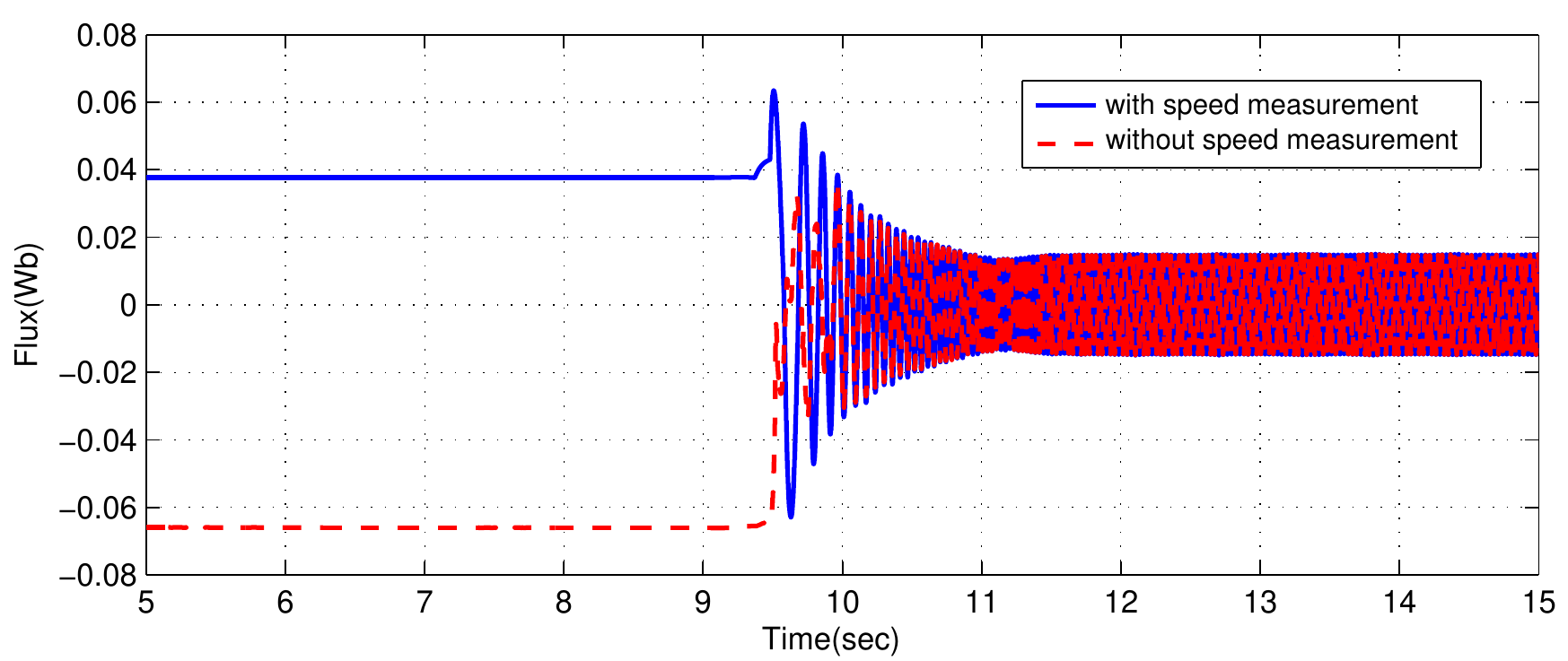}
	\label{psi_ra_est_exp_im}}

\subfloat[$\psi_{r \beta_s}$]{	\includegraphics[width = 0.85\linewidth]{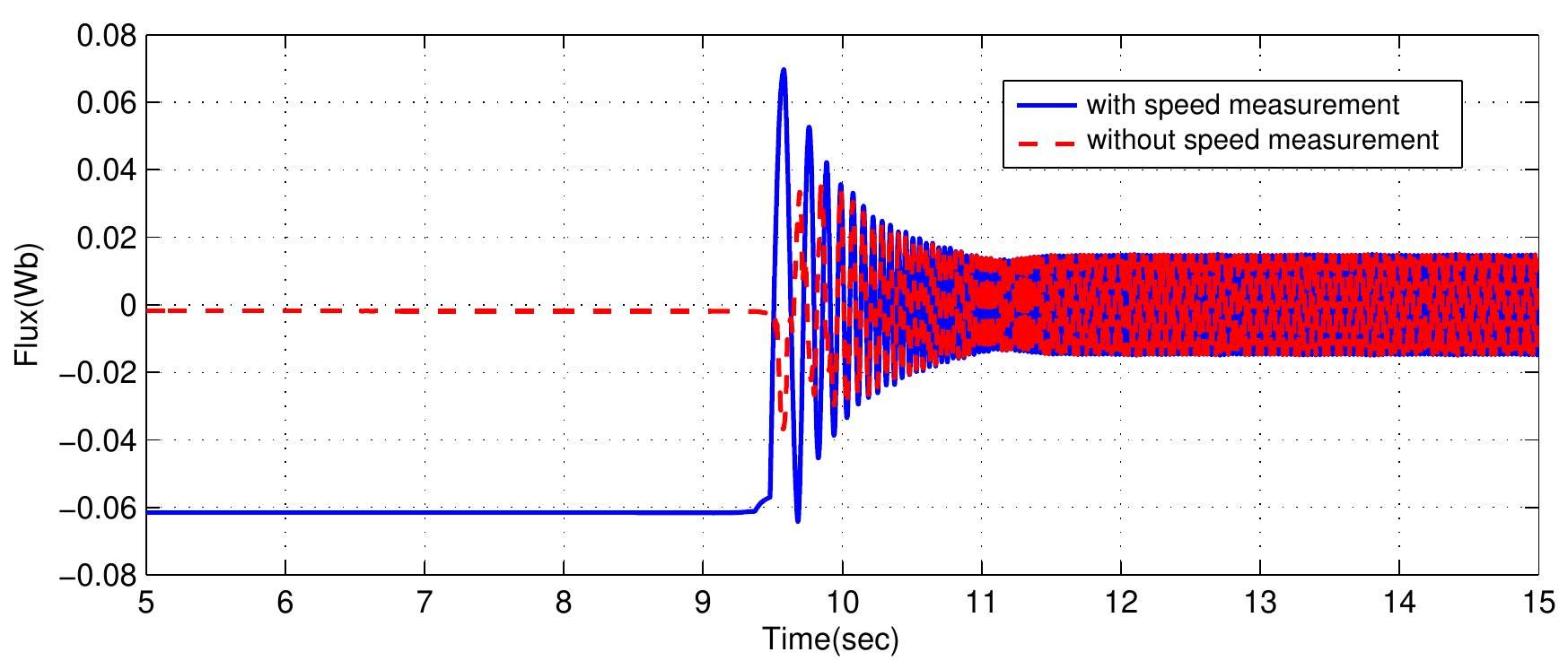}
	\label{psi_rb_est_exp_im}}
    \caption{Rotor flux estimation}
\end{figure}

\begin{figure}[!hb]
	\centering
	\includegraphics[width = 0.85\linewidth]{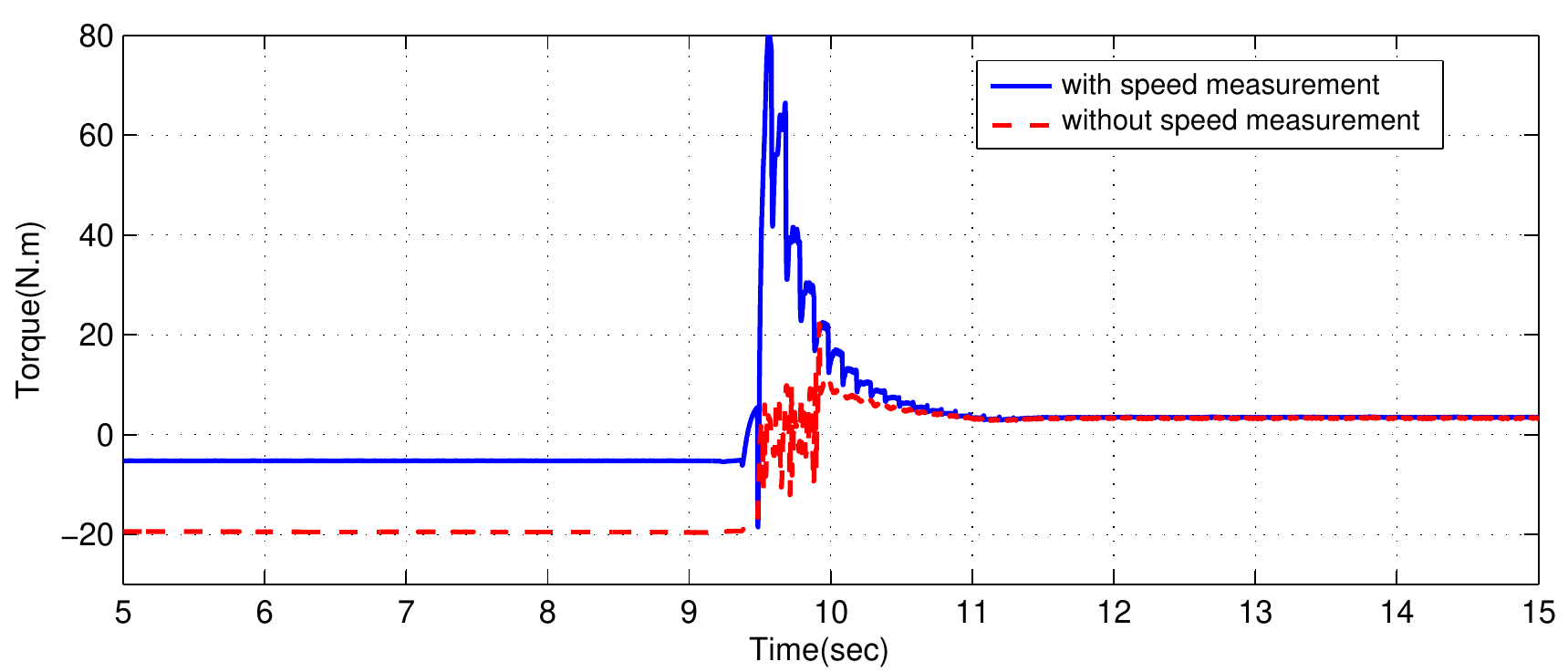}
	\caption{Load torque estimation}
	\label{cr_esti_exp_im}
\end{figure}

In the absence of torque and flux measurement, it will considered that the observer having speed measurement as input converges to the correct value in steady-state. The estimated rotor fluxes are shown on the figures \ref{psi_ra_est_exp_im} for $\hat{\psi}_{r \alpha}$ and \ref{psi_rb_est_exp_im} for $\hat{\psi}_{r\beta}$: at standstill ($\omega_e = 0$), with zero stator frequency ($\omega_s = 0$), the estimated flux without speed measurement is different from that estimated with speed measurement. The same observation is noticed for the load torque, shown on the figure \ref{cr_esti_exp_im}, and the rotor speed shown on the figure \ref{omega_est_exp_im}. However, for non-zero speed and stator frequency, the estimated states without speed measurement converge to those with speed measurement. Therefore, both simulation and experimental results are in coherence with the observability study of Section 4.



\section{Conclusions}
The local observability of electric drives in view of sensorless control has been studied throughout this paper. The study shows that the observability is not guaranteed in some operating conditions, such as standstill operation for the synchronous machines, and low-frequency input voltages for the induction machines.

One interesting feature of the observability approach applied in this paper is the physically insightful analytic formulation of observability conditions, thanks to the rank criterion. These conditions can be easily interpreted and tested in real time. Furthermore, the study resulted in the definition of two concepts: the \emph{observability vector} for the SMs and the \emph{unobservability line} for the IMs. The validity of the observability conditions is confirmed by numerical simulations and experimental data, using an extended Kalman observer.

Finally, the following conclusions on electric drives observability have been drawn:
\begin{itemize}
\item The observability of DC machines is guaranteed if the field flux is not zero.
\item The observability of SMs is guaranteed if the angular velocity of the \emph{observability vector} in the rotor reference frame is different from the angular velocity of the rotor.
\item The observability of IM is guaranteed if the stator voltage frequency is not zero.
\end{itemize}

Several solutions are proposed in the literature to recover the observability in the critical situations, based on the existing observability conditions \cite{zgorski_sled_12,gaeten_sled_15}. We believe that the new results presented in this paper are so useful for engineers and researchers seeking better understanding of sensorless drives performance. For instance, based on the observability analysis of the wound rotor synchronous machine, we proposed a technique that consists of superimposing a high-frequency current with the rotor current at standstill in order to improve the observer performance. Other techniques can be investigated for different machines.

\appendix

\section{Machines parameters}
\label{appendix_a}
The WRSM parameters are the following:
\begin{eqnarray}
&& R_s =  0.01~\Omega~~;~~ L_d = 0.8~mH ~~~;~~ L_q = 0.7~mH;\\
&& R_f = 6.5~\Omega~~~;~~ L_f =850~mH~~;~~ M_f = 5.7~mH \\
&& J = 10^{-2} kg.m^2~~;~~ p=2
\end{eqnarray}
The IM parameters are the following:
\begin{eqnarray}
&& R_s = 2.8~m\Omega ~~~~~~~~~~~~;~~L_s = 9.865\times10^{-5} ~H  \\
&& R_r = 1.5~m\Omega ~~~~~~~~~~~~;~~ L_r = 1.033\times 10^{-4}~H \\
&& M = 9.395\times10^{-5}~H ~~;~~ p = 4\\
&& J = 10^{-2} kg.m^{2}~~~~~~~~~~;~~f_v = 10^{-4}~N.m.s.rad^{-1}
\end{eqnarray}

\bibliographystyle{apalike}
\bibliography{bib_drives_observability}

\end{document}